\documentclass{article}

\usepackage{amsmath}
\usepackage{amsthm}
\usepackage{amssymb}
\usepackage{graphicx}
\usepackage{cellspace}
\usepackage{geometry}
\usepackage{pgf,tikz,pgfplots}
\usepackage[all]{xy}
\usepackage{mathrsfs}
\usetikzlibrary{arrows}
\usepackage[utf8]{inputenc}
\usepackage{titlesec}
\usepackage{multicol}
\usepackage{algorithm}
\usepackage{algpseudocode}
\usepackage[normalem]{ulem}
\usepackage[T1]{fontenc} 
\usepackage[utf8]{inputenc} 
\usepackage{charter}
\usepackage{array}
\usepackage{calc}
\usepackage{tikz-3dplot}
\usepackage{wasysym}
\usepackage{marvosym}
\usetikzlibrary[patterns]
\usepackage{mdframed} 
\pgfplotsset{compat=1.17}

\definecolor{rouge}{RGB}{204, 15, 39}
\definecolor{bleu}{RGB}{60, 86, 164}
\definecolor{vert}{RGB}{10, 195, 118}
\definecolor{orange}{RGB}{245, 85, 63}
\definecolor{gris}{RGB}{69, 94, 93}
\definecolor{noir}{RGB}{0, 0, 0}
\definecolor{jaune}{RGB}{255, 195, 42}

\begingroup
\makeatletter
\@for\theoremstyle:=definition,remark,plain\do{%
	\expandafter\g@addto@macro\csname th@\theoremstyle\endcsname{%
		\addtolength\thm@preskip\parskip
	}%
}
\endgroup

\newtheorem*{thm*}{Theorem}
\newtheorem*{pties*}{Properties}
\newtheorem*{prop*}{Proposition}
\newtheorem*{cor*}{Corollary}
\newtheorem*{defn*}{Definition}

\newtheorem{thm}{Theorem}[section]
\newtheorem{lem}[thm]{Lemma}

\newtheorem{prop}[thm]{Proposition}

\newtheorem{cor}[thm]{Corollary}
\newtheorem{defn}[thm]{Definition}

\theoremstyle{remark}
\newtheorem{rem}[thm]{Remark}

\newtheorem{ex2}[thm]{Example}
\newenvironment{ex}{\begin{ex2}}{\hfill $\diamondsuit$ \end{ex2}}

\newtheorem*{pro}{\textit{Proof}}
\newtheorem*{proth}{\textit{Proof of the theorem}}
\newtheorem*{propr}{\textit{Proof of the proposition}}
\newtheorem{nota}[thm]{Notation}
\newenvironment{prr}{\begin{pro}}{\hfill $\square$ \end{pro}}

\newcommand{\bC}{\mathbb{C}}

\newcommand{\bP}{\mathbb{P}}
\newcommand{\bQ}{\mathbb{Q}}
\newcommand{\bR}{\mathbb{R}}

\newcommand{\bZ}{\mathbb{Z}}

\newcommand{\cE}{\mathcal{E}}
\newcommand{\cF}{\mathcal{F}}
\newcommand{\cG}{\mathcal{G}}
\newcommand{\cH}{\mathcal{H}}

\newcommand{\cL}{\mathcal{L}}

\newcommand{\cO}{\mathcal{O}}

\newcommand{\cT}{\mathcal{T}}

\newcommand{\sB}{\mathscr{B}}

\newcommand{\sE}{\mathscr{E}}

\newcommand{\im}{\textup{Im}}

\newcommand{\di}{\textup{div}}

\newcommand{\rk}{\textup{rk}}

\newcommand{\Sp}{\textup{Spec}}

\newcommand{\Sym}{\textup{Sym}}

\newcommand{\Span}{\textup{Span}}

\usepackage{booktabs}
\usepackage{setspace}
\definecolor{link}{RGB}{71,81,114}
\usepackage{hyperref}                
\hypersetup{
	colorlinks = true,
	allcolors=link}

\title{A combinatorial description of stability for toric vector bundles}
\author{Lucie DEVEY}
\date{}
\newcommand{\Addresses}{{
		\bigskip
		\footnotesize\textsc{Institut für Mathematik, Goethe–Universität Frankfurt, 60325 Frankfurt am Main, Germany\\
			Institut Fourier, Université Grenoble Alpes, 38400 Saint Martin d’Hères, France}\par\nopagebreak
		\textit{E-mail address}, \texttt{lucie.devey@univ-grenoble-alpes.fr}
}}

\begin{document}

\maketitle

\begin{abstract}
	The aim of this paper is to discuss a combinatorial characterisation of stability for toric vector bundles (or equivariant reflexive sheaves) in the terms of their parliaments of polytopes, a generalization of moment polytopes for toric vector bundles by Di Rocco, Jabbusch and Smith. We also define subparliaments of polytopes and identify them with parliaments of equivariant subbundles.
\end{abstract}

\section*{Introduction}

Stability has been introduced in the perspective of classifying vector bundles, by constructing moduli spaces of semistable sheaves. Using the Harder--Narasimhan filtration, any vector bundle $\cE$ can be built up from semistable coherent sheaves. In 1963, Mumford extended the works of Grothendieck (see \cite{G} Theorem 2.1) and Atiyah (see \cite{A57}) by classifying vector bundles on curves (see \cite{M63}). In this paper, we consider a generalisation of Mumford-stability to vector bundles on varieties of any dimension, called slope-stability and introduced in 1972 by Takemoto (\cite{T72}). The problem of classification of vector bundles is still completely open. The aim of this paper is to study the toric case.

Consider a toric variety $X$ with torus $T$. A toric vector bundle $\cE$ over $X$ is a locally free $\cO_X$-module of finite rank $r$, equipped with a $T$-action such that the projection map $\pi:\Sp(\Sym\,\cE)\to X$ is $T$-equivariant and $T$ acts linearly on the fibres. We denote the fiber over the identity by $E$. In \cite{DrJS}, Sandra Di Rocco, Kelly Jabbusch and Gregory G. Smith generalized the construction of the moment polytope of a line bundle, associating to any equivariant vector
bundle $\cE$ its parliament of polytopes $PP_\cE$ : a collection of convex polytopes $(P_e)_{e\in G(\cE)}$ indexed by the elements in the ground set of a matroid $M(\cE)$ associated to $\cE$ representable in the vector space $E\simeq\bC^r$. The wealth of information about a toric vector bundle $\cE$ contained in its parliament is astounding. For instance, lattice points in the parliament of polytopes for $\cE$ correspond to a torus-equivariant generating set for the space of global sections of $\cE$ (see \cite{DrJS} Proposition 1.1). In addition, we can also recover some positivity properties of $\cE$, such as global generatedness (see \cite{DrJS} Theorem 1.2), ampleness (see \cite{DrJS} Corollary 6.7), or bigness (see \cite{N} Theorem 7.5). 

Our main result is a first step to having an algorithm (as wished in \cite{DDK1}) for the slope-stability of toric vector bundles. We start by defining in Definition \ref{average} for any equivariant saturated subsheaf $\cF$ of a toric vector bundle $\cE$, a polytope called the average polytope $P_{\cF}$ which can be visualized on the parliament of $\cE$. It is the moment polytope of $c1(\cF)/\rk(\cF)$. In Definition \ref{order}, for any polarisation $\alpha$ of $X$, we define a total order $<_\alpha$ between polytopes, such that comparing average polytopes $P_{\cF_1}$ and $P_{\cF_2}$ corresponds to comparing the $\alpha$-slopes of their respective sheaves $\cF_1$ and $\cF_2$. To check stability, we need to compare the slope of $\cE$ with the slope of its subsheaves. It has been proved, first that equivariant saturated subsheaves of $\cE$ are in one-to-one correspondence with vector subspaces $F\subset E$ and second that it is enough to check the slopes of those subsheaves are less than the slope of $\cE$ (see \cite{DDK1} or \cite{HNS}). These nice results allowed the same authors to have a finite check for the stability of tangent bundles. In this paper, we give the existence of a finite matroid $M(\cE)^S$ on which the flats correspond to a finite family of subsheaves $\cF$ sufficiently varied to check stability of any toric vector bundle.
\begin{thm*}[\ref{thm:existence}] 
	Let $X$ be a smooth complete toric variety. For any toric vector bundle $\cE$ on $X$, there exists a finite matroid $M(\cE)^S$ such that $\cE$ is $\alpha$-(semi)stable if and only if, for any nonzero flat $f\subsetneq G(\cE)^S$ of $M(\cE)^S$, we have
	$$P_\cF <_\alpha P_\cE\quad \text{(resp. }P_\cF \leq_\alpha P_\cE)\ ,$$
	where $\cF$ is the equivariant saturated sheaf corresponding by Theorem 2.4 to the linear subspace $\Span(f)\subset E$.
\end{thm*}

Moreover, we give a finite check of stability for any toric bundle of rank less or equal to $3$: in that case, the matroid of any parliament of polytopes of $\cE$ checks stability. We caution the reader that the first version of the article contains a mistake, we claimed that this same finite check is working for any rank, it is false.

In another part of the article, we prove that any equivariant subbundle $\cF$ of a toric vector bundle $\cE$ corresponds to a flat of a matroid $M(\cE)$ compatible with the Klyachko filtrations of $\cE$. We define subparliaments of polytopes and identify them with parliaments of equivariant subbundles.
\begin{cor*}[\ref{subp}]
	The subparliaments of the parliaments of $\cE$ are the parliaments of the equivariant subbundles of $\cE$.
\end{cor*}

The last part of this article consists of exploring the geometric information that can
be reconstructed from the parliament of polytopes of a toric vector bundle. In particular, we translate results of Payne (see \cite{P} Proposition 3.4, Corollary 3.6) and Klyachko (see \cite{K} Corollary 1.2.4) in terms of parliaments of polytopes.

\begin{prop*}[\ref{isom2}]
	The data of the parliament of polytopes of a globally generated equivariant vector bundle $\cE$, up to translation of each direct component and quotiented by $GL_r(\bC)$, enables us to reconstruct the isomorphism class of the vector bundle $\cE$.
\end{prop*} 

We start in Section \ref{s:recall} by recalling facts about matroids, giving the construction of
parliaments of polytopes from \cite{DrJS} and fixing notation.

Section \ref{s:stability} furnishes the definition of average polytopes that allows us to visualize the
slope of an equivariant saturated sheaf. It also contains the most important result of
the article (Theorem \ref{thm:existence}): the existence, for any equivariant vector bundle $\cE$, of a finite matroid checking combinatorially the slope-stability of $\cE$ using one of its parliaments of polytopes. The cases of tangent bundles and of vector bundles of rank $3$ are also discussed in this section.

In Section \ref{s:sub}, given the parliament of polytopes of a toric vector bundle, we describe
the parliaments of its equivariant subbundles.

In Section \ref{s:rest}, we treat the stability of the restriction of an equivariant bundle $\cE$ to
a torus invariant curve in terms of the parliament $PP_\cE$. We also give examples of $\alpha$-(semi)stable equivariant vector bundle with non (semi)stable restrictions to torus invariant curves.

In Section \ref{s:A}, we discuss the definition of parliament of polytopes and state what data
is encoded in the parliament of an equivariant bundle.

\section*{Acknowledgements}

I would like to express my deep gratitude to my thesis advisors Alex Küronya and Catriona Maclean. I warmly thank Klaus Altmann, Sebastien Boucksom, Christian Haase, Milena Hering, Bivas Khan, Diane Maclagan, Chris Manon, Joaquim Roé, Karin Schaller, Greg Smith and Martin Ulirsch.
This work has been partially supported by the LabEx PERSYVAL-Lab (ANR-11-LABX-0025-01) funded by the French program Investissement d'avenir and by the Deutsche Forschungsgemeinschaft (DFG, German Research Foundation) TRR 326 Geometry and Arithmetic of Uniformized Structures, project number 444845124.

\setlength\parindent{0pt}

\section{Parliaments of polytopes and matroid terminology}\label{s:recall}

Consider a smooth complete toric variety $X$ of dimension $d$ with fan $\Sigma$ and torus $T$. The lattice of characters and its dual lattice are denoted by $M$ and $N$. Some matroids thereafter will have the same notation $M$, we hope that it will not cause confusions. We denote by $\Sigma(k)$ the cones of dimension $k$. Let $\cE$ be a given rank-$r$ equivariant vector bundle over $X$ with fiber over the identity of the torus $E$. Let us denote by $\{v_0,\hdots,v_{n-1}\}$ the set of vectors generating the rays $\rho_i$ of $\Sigma$. (See Chapter \cite{CLS} for more fundamentals on toric varieties.) Parliaments of polytopes were introduced in \cite{DrJS} in order to give explicit polyhedral interpretations of properties to equivariant vector bundles. The parliaments of polytopes $PP_\cE$ of $\cE$ are composed of at least $r$ polytopes linked with some combinatorial data: the polytopes are labelled by elements of a representable matroid of rank $r$.

We first recall some matroid terminology. See \cite{K14} for more fundamentals on matroids.

\subsection{Matroid terminology}

Matroids are a generalization of the notion of linear independence in vector spaces. 

\begin{defn}
	A matroid $M$ is the data of a finite set $G$ and a collection $\sB$, of subsets of $G$, called bases, satisfying the following properties: 
	\begin{itemize}
		\item [B1:] $\sB$ is nonempty ;
		\item [B2:] (basis exchange property) If $A,B\in\sB$ are distinct and $a\in A\setminus B$, then there exists $b\in B\setminus A$ such that $(A\setminus \{a\})\cup \{b\}\in \sB$ .
	\end{itemize}
	We call $G$ the ground set and $\sB$ the set of bases of the matroid $M=(G,\sB)$.
\end{defn}

\begin{defn}
	An isomorphism of matroids $\varphi:M_1=(G_1,\sB_1)\to M_2=(G_2,\sB_2)$ is a bijection from $G_1$ to $G_2$ such that $$A\in\sB_1\Leftrightarrow \varphi(A)\in\sB_2\ .$$
\end{defn}

\begin{defn}
	A representable matroid of rank $r$ is a matroid isomorphic to $$M=(G,\sB)\ ,$$ where $G$ is a finite subset of some dimension $r$-vector space $E$ and $\sB$ is the set of bases of $E$ formed by vectors in $G$. We say that $M$ is represented in $E$.
\end{defn}

\begin{defn}
	A flat of a matroid $(G,\sB)$ represented in $E$, is a subset $f\subset G\subset E$ such that $$\Span(f)\cap G=f\ ,$$ where $\Span(f)\subset E$ is the subspace spanned by the vectors in $f$.
\end{defn}

\subsection{Parliaments of polytopes}\label{ss:parliament}

The construction of parliaments of polytopes of $\cE$ requires the Klyachko classification.

\begin{thm}[Theorem 0.1.1 of \cite{K}]\label{K} There is an equivalence of categories between the category of rank-$r$ equivariant vector bundles $\cE$ on $X$ and the category of $n$ compatible decreasing $\bZ$-filtrations $(E^i(j))_{\rho_i\in\Sigma(1)}$
	$$E^i(j)=\left\{\begin{array}{ll}
		E &\text{if }j\leq A^i_{1}\\
		H_i &\text{if }A^i_{1}<j\leq A^i_{2}\\
		\ldots &\\
		\Span(u_i)  &\text{if }A^i_{r-1}<j\leq A^i_{r}\\
		\{0\} &\text{if } A^i_{r}<j
	\end{array}\right.\,.$$
	of a dimension $r$ $\bC$-vector space $E$ together with a compatibility condition as follows. There exist decompositions of $E$ into $1$-dimensional vector spaces $L_u^\sigma$ satisfying
	
	\begin{equation}\tag{CC}\label{CC}
		\forall\,\sigma\in\Sigma(d),\ \exists\, (L_u^\sigma)_{u\in\textbf{u}(\sigma)} \text{ s.t. }E=\bigoplus_{u\in\textbf u(\sigma)}L_u^\sigma\text{ and }\ \forall\,\rho_i\preceq
		\sigma, \; E^i(j)=\sum_{\langle u,v_i\rangle \geq j}L_u^\sigma\ .
	\end{equation} 
\end{thm}

\begin{rem}\label{eqcat}
	There is an equivalence of categories between the category of equivariant reflexive sheaves $\cF$ on $X$ of rank $r$, and the category of $n$ compatible decreasing $\bZ$-filtrations of a $r$-dimensional $\bC$-vector space $F$ where we do not impose the compatible condition (see Theorem 5.19 \cite{P01}).
\end{rem}

\begin{defn}
	The $r$ points in $\textbf{u}(\sigma)$ are called the associated characters of $\sigma$.
\end{defn}

\begin{rem}
	The compatibility condition (\ref{CC}) implies the existence, for every maximal cone $\sigma\in\Sigma(d)$, of a basis $B_\sigma$, given by a generator of $L_u^\sigma$ for each $u\in\textup{u}(\sigma)$. We call $B_\sigma$ a compatible basis. Nevertheless, given a maximal cone, a compatible basis may not be unique (see Example 4.4 of \cite{DrJS}) although $\textbf{u}(\sigma)$ always is. 
\end{rem}

\begin{rem}
	An important feature is that for any maximal cone $\sigma\in\Sigma(d)$, the equivariant bundle $\cE$ splits equivariantly on $U_\sigma$ as 
	$$\cE|_{U_\sigma}\simeq\bigoplus_{u\in\textbf{u}(\sigma)}\cO_X\left(\di (u)\right)|_{U_\sigma}\ .$$
\end{rem}

\noindent

A parliament of polytopes of a vector bundle $\cE$ is a certain collection of polytopes
$$\left\{(P_e,e)\ \middle|\ e\in G(\cE)\subseteq E\right\}\,,$$ where $G(\cE)$ is the ground set of a matroid $M(\cE)$ defined in the following manner.
The matroids, called $M(\cE)$, associated to $\cE$ can be seen as minimal matroids generating
$$L(\cE)=\left\{\bigcap_{i\in\{0,\hdots,n\},}E^i(j_i)\,\middle|\;(j_i)_{i\in\{0,\hdots, n\}}\in\bZ^{n+1}\right\}$$
as a meet-subsemilattice (as a partially ordered subset which has a meet i.e. a greatest lower bound). Despite Proposition 3.1 \cite{DrJS}, Algorithm 3.2 does not produce a unique matroid up to isomorphism (i.e. the procedure is ambiguous in some cases)\footnote{see Examples \ref{cex} and \ref{cex2} for counterexamples}.
However, the ground sets $G(\cE)$ generating the matroids $M(\cE)$ are the possible outputs of Algorithm 3.2 of \cite{DrJS}.
{\begin{algorithm}
		\renewcommand{\thealgorithm}{3.2}
		\label{3.2}
		\begin{algorithmic}
			\State $r\gets$ the dimension of the largest linear subspace of $L(\cE)$
			\State $G\gets\emptyset$
			\For{$k=1$ to $r$}
			\For{$k$-dimensional linear subspace $V\in L(\cE)$}
			\State $G'\gets G\cap V$
			\If{$\textup{Span}(G')\subsetneq V$} 
			\State $G\gets G\sqcup CB_V$, $\quad CB_V$ is a basis of a complement to $\textup{Span}(G')$ in $V$ 			\EndIf
			\EndFor
			\EndFor
			\State \textbf{return} $G$
		\end{algorithmic}
		\caption{\quad Computing $G(\cE)$}
\end{algorithm}}

Additionally, if $G_E^1$ and $G_E^2$ are ground sets resulting from Algorithm 3.2 applied to $L(\cE)$ then the same number of elements is added at Step $V$. Indeed, $CB_V^1,CB_V^2$ are bases of a complement to $\textup{Span}(G'^1),\textup{Span}(G'^2)$, which have same dimension 
$$\dim(\textup{Span}(G'^\varepsilon))=\dim\left(\textup{Span}\left(V'\cap V\right)_{V'\cap V\subsetneq V}\right)\quad \forall\,\varepsilon\in\{1,2\}\,$$
in $V$. The equality holds because, as $L(\cE)$ is stable by intersection, the vector spaces $V'\cap V$ that are different from $V$, are steps of the algorithm and are already generated by $G'^\varepsilon$ at Step $V$. Therefore, there exists a map $\varphi$ sending $G_E^1$ to $G_E^2$, we call it of type $(\star)$.

\begin{defn}\label{iso}
	A map $\varphi:M(G_E^1)\to M(G_E^2)$ of type $(\star)$ is a bijection 
	\begin{equation*}
		\varphi:G_E^1\to G_E^2
	\end{equation*}
	which respects Algorithm 3.2 for $\cE$, that is to say which sends $CB_V^1$ to $CB_V^2$ at Step $V$ of Algorithm 3.2.
\end{defn} 

We use the following notation for a ground set $G(\cE)$ of $M(\cE)$ 
\begin{center}
	$G(\cE)=\{e_0,\hdots,e_l\}\quad$ for some $l\geq r$\ . 
\end{center}

\begin{defn}\label{def}
	Any output $G(\sE)$ from Algorithm 3.2 gives rise to a parliament of polytopes of $\cE$. These are the sets of indexed polytopes $(P_e)_{e\in G(\cE)}$ defined as 
	$$P_e:=\left\{m\in M_\bR\,\middle|\;\forall\rho_i\in\Sigma(1),\:\langle m,v_i\rangle \leq\max\left\{j\,\middle|\;e\in E^i(j)\right\}\right\}\ ,$$
	where $G(\cE)$ is seen modulo isomorphism of matroids of type $(\star)$.
\end{defn}

\begin{ex}\label{cex}[Communicated by Diane Maclagan]
	To a fixed toric vector bundle, we associated a finite number of parliaments of polytopes, corresponding to each matroid $M(\cE)$ (modulo isomorphism) obtained by Algorithm 3.2. Consider the canonical basis $\{e_1,e_2,e_3\}$ of $\bC^3$.
	
	Let $X$ be the projective plane $\bP^2$ and let $\cE$ be the toric bundle of rank $3$ defined by the following filtrations 
	$$E^i(j)=\left\{\begin{array}{ll}
		\bC^3 &\text{if }j\leq 0\\
		V_i &\text{if }0<j\leq 1\\
		\{0\} &\text{if } 1<j
	\end{array}\right.\,,$$
	where $V_0=\Span(e_1)$, $V_1=\Span(e_1+e_2)$ and $V_2=\Span(e_2,e_3)$.
	The meet-subsemilattice is $$L_1(\cE)=\left\{\{0\},\Span(e_1),\Span(e_1+e_2),\Span(e_2,e_3),\bC^3\}\right\}\ .$$  Two non-isomorphic matroids are given by $G_1=\{e_1,e_1+e_2,e_2,e_3\}$ and $G_2=\{e_1,e_1+e_2,e_2+e_3,e_3\}$.
	Here the non-uniqueness comes from $L_1(\cE)$ not being stable by sums. 
\end{ex}
	
\begin{ex}\label{cex2}
	A second example can be found on $\bP^2$, let $\cE$ be defined by the following filtrations 
	$$E^i(j)=\left\{\begin{array}{ll}
		\bC^3 &\text{if }j\leq 0\\
		V_i &\text{if }0<j\leq 1\\
		\Span(e_1) &\text{if }1<j\leq 2\\
		\{0\} &\text{if } 2<j
	\end{array}\right.\,,$$
	where $V_0=\Span(e_1,e_2)$, $V_1=\Span(e_1,e_3)$ and $V_2=\Span(e_1,e_2+e_3)$.
	The meet-subsemilattice is 
	$$L_2(\cE)=\left\{\{0\},\Span(e_1),\Span(e_1,e_2),\Span(e_1,e_3),\Span(e_1,e_2+e_3),\bC^3\}\right\}\ .$$ Two non-isomorphic matroids are given by $G_1=\{e_1,e_2,e_2+e_3,e_3\}$ and $G_2=\{e_1,e_1+e_2,e_2+e_3,e_3\}$. 
	Here the non-uniqueness comes from $L_2(\cE)\cup \{G_2\}$ not being stable by iterated sums and intersections. 
\end{ex}

For simplicity, we extend the definition of parliaments of toric bundles (Definition \ref{def}) to preparliaments of equivariant reflexive sheaves\footnote{A coherent sheaf $\cF$ on $X$ is reflexive if $\cF$ is isomorphic to its double dual $\cF^{**}$.}. Perling, in \cite{P01} Theorem 5.19, extends the Klyachko classification as follows.

\begin{prop}\label{perl}
	There is an equivalence of categories between the category of rank-$r$ equivariant reflexive sheaves $\cE$ on $X$ and the category of $n+1$ decreasing $\bZ$-filtrations $(E^i(j))_{\rho_i\in\Sigma(1)}$ of a dimension $r$-$\bC$-vector space $E$
	$$E^i(j)=\left\{\begin{array}{ll}
		E &\text{if }j\leq A^i_{1}\\
		H_i &\text{if }A^i_{1}<j\leq A^i_{2}\\
		\ldots &\\
		\Span(u_i)  &\text{if }A^i_{r-1}<j\leq A^i_{r}\\
		\{0\} &\text{if } A^i_{r}<j
	\end{array}\right.\,.$$
\end{prop}

We define accordingly the notion of preparliaments of polytopes.

\begin{defn}\label{ppp}
	Let $\cE$ be an equivariant reflexive sheaf on $X$.
	A preparliament of polytopes of $\cE$ is composed of a set of indexed polytopes $(P_g)_{g\in G(\cE)}$ defined as 
	$$P_g:=\left\{m\in \bR^d\,\middle|\;\forall\rho_i\in\Sigma(1),\:\langle m,v_i\rangle \leq\max\left\{j\,\middle|\;g\in E^i(j)\right\}\right\}\ ,$$
	where the ground set $G(\cE)$ is obtained by Algorithm 3.2 and is seen modulo isomorphism of matroids of type $(\star)$.
\end{defn}

A preparliament that satisfies the compatibility condition (\ref{CC}) is a parliament, and the corresponding equivariant reflexive sheaf is a toric bundle. We keep the notation $PP_\cE$ for preparliaments.

\begin{ex}\label{tp2}
	Consider $X=\bP^2$ and its tangent bundle $\cE=\cT_X$. The $\bZ$-filtrations are $$\text{for }i=0,1,2,\;E^i(j)=\left\{\begin{array}{ll}
		\bC^2 &\text{if }j<0\\
		\Span(v_i)  &\text{if }0\leq j<1\\
		\{0\} &\text{otherwise}
	\end{array}\right..$$
	The ground set is $G(\cE)=\{v_0,v_1,v_2\}=\Sigma(1)$ and the parliament $PP_\cE$ is the following.
	
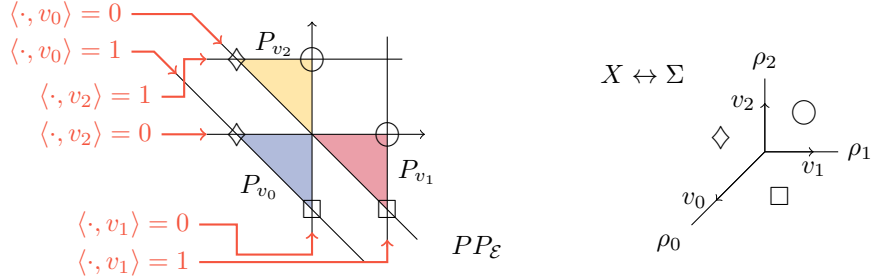
\vspace{-5pt}
	$$\begin{tikzpicture}[scale=1]
		\fill [jaune!40] (-1,1) -- (0,1) -- (0,0);
		\fill [bleu!40] (-1,0)-- (0,-1) -- (0,0);
		\fill [rouge!40] (1,-1) -- (1,0) -- (0,0);
		\draw (-1.4,1) -- (1.2,1);
		\draw [->] (-1.4,0) -- (1.5,0);
		\draw (1,-1.5) -- (1,1.3);
		\draw [->] (0,-1.5) -- (0,1.5);
		\draw (-1.4,1.4) -- (1.4,-1.4);
		\draw (-2.1,1.1) -- (.7,-1.7);   
		\draw [->,orange,line width=.8] (-2.4,1.6) node[left] {$\langle\cdot,v_0\rangle=0$} -- (-1.6,1.6) -- (-1.2,1.2);
		\draw [->,orange,line width=.8] (-2.4,1.1) node[left] {$\langle\cdot,v_0\rangle=1$} -- (-2.1,1.1) -- (-1.8,.8);
		\draw [->,orange,line width=.8] (-2,.5) node[left] {$\langle\cdot,v_2\rangle=1$} -- (-1.7,.5) -- (-1.7,1) -- (-1.4,1);
		\draw [->,orange,line width=.8] (-2,0)  node[left] {$\langle\cdot,v_2\rangle=0$} -- (-1.4,0);
		\draw [->,orange,line width=.8] (-1.5,-1.7) node[left] {$\langle\cdot,v_1\rangle=1$} -- (1,-1.7)  -- (1,-1.4);
		\draw [->,orange,line width=.8] (-1.5,-1.2) node[left] {$\langle\cdot,v_1\rangle=0$} -- (-1,-1.2) -- (-1,-1.6) -- (0,-1.6) -- (0,-1.3);
		\draw (-.5,1.2) node {$P_{v_2}$};   
		\draw (1.4,-.5) node {$P_{v_1}$};
		\draw (-.7,-.7) node {$P_{v_0}$};
		\draw (2.2,-1.5) node {$PP_\cE$};
		\draw (-1,1) node {$\diamondsuit$};
		\draw (-1,0) node {$\diamondsuit$};
		\draw (0,-1) node {$\square$};
		\draw (1,-1) node {$\square$};
		\draw (1,0) node {$\bigcirc$};
		\draw (0,1) node {$\bigcirc$};
	\end{tikzpicture}\quad \qquad\begin{tikzpicture}[scale=.65]
		\draw (-2.5,1.6) node {$X\leftrightarrow\Sigma$};
		\draw (0,0) -- (1.5,0) node[right] {$\rho_1$};
		\draw (0,0) -- (0,1.5)node[above] {$\rho_2$};
		\draw (-1.5,-1.5) node [below left]  {$\rho_0$}-- (0,0);
		\draw [->] (0,0) -- (1,0) node[below] {$v_1$};
		\draw [->] (0,0) -- (0,1) node[left] {$v_2$};
		\draw [<-] (-1,-1) node [left]  {$v_0$} -- (0,0);
		\draw (-.9,.3) node {$\diamondsuit$};
		\draw (.3,-.9) node {$\square$};
		\draw (.8,.8) node {$\bigcirc$};
		\draw [white] (0,-2) -- (0,-2.7);
	\end{tikzpicture}$$
	We associate to each maximal cone $\sigma$ a symbol (say $\square$) and we represent the associated characters $\textbf u(\sigma)$ by $r$ symbols $\square$ on the parliament. The compatibility condition (\ref{CC}) is verified.
\end{ex}

\begin{ex}
	A contrario, the equivariant reflexive sheaf on $\bP^3$ defined by the filtrations
	$$E^i(j)=\left\{\begin{array}{ll}
		\bC^3 &\text{if }j\leq 0\\
		\textup{Span}(u_i,u) &\text{if }0<j\leq 1\\
		\textup{Span}(u_i)  &\text{if }1<j\leq 2\\
		\{0\} &\text{if } 2<j
	\end{array}\right.\,,$$
	where $u_0,u_1,u_2,u_3$ belongs to a plan $P$ and $u$ does not,
	$$
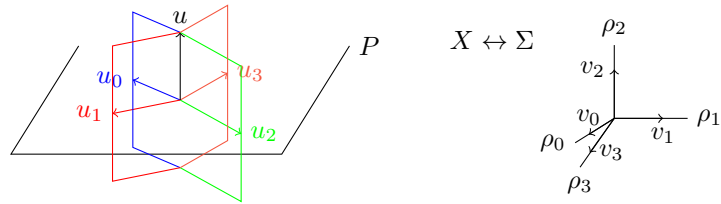
\begin{tikzpicture}[scale=.9]
		\draw (0,1.8) -- (-1,0.2) -- (3,0.2) -- (4,1.8) node [right] {$P$};
		\draw [->, color=blue] (1.5,1) -- (.8,1.3) node [left]{$u_0$};
		\draw [->,color=red] (1.5,1) -- (.5,.8) node[left] {$u_1$};
		\draw [->,color=green] (1.5,1) -- (2.4,.5) node[right] {$u_2$};
		\draw [->,color=orange] (1.5,1) -- (2.2,1.4) node[right] {$u_3$};
		\draw [->] (1.5,1) -- (1.5,2) node[above] {$u$};
		\draw [blue] (1.5,2) -- (.8,2.3) -- (.8,.3) -- (1.5,0);
		\draw [red] (1.5,2) -- (.5,1.8) -- (.5,-.2) -- (1.5,0);
		\draw [green] (1.5,2) -- (2.4,1.5) -- (2.4,-.5) -- (1.5,0);
		\draw [orange] (1.5,2) -- (2.2,2.4) -- (2.2,.4) -- (1.5,0);
	\end{tikzpicture}\qquad
	\begin{tikzpicture}[scale=.65]
		\draw (-2.5,1.6) node {$X\leftrightarrow\Sigma$};
		\draw (0,0) -- (1.5,0) node[right] {$\rho_1$};
		\draw (0,0) -- (0,1.5) node[above] {$\rho_2$};
		\draw (0,0) -- (-.7,-1) node[below] {$\rho_3$};
		\draw (-.8,-.5) node [left]  {$\rho_0$}-- (0,0);
		\draw [->] (0,0) -- (1,0) node[below] {$v_1$};
		\draw [->] (0,0) -- (0,1) node[left] {$v_2$};
		\draw [->] (0,0) -- (-.5,-.7) node[right] {$v_3$};
		\draw [<-] (-.52,-.33) node [above]  {$v_0$} -- (0,0);
		\draw [white] (0,-1.7) -- (0,-1.5);
	\end{tikzpicture}
	$$
	is not a toric bundle, but we can still construct its preparliament of polytopes. It contains five polytopes indexed by $u_0$, $u_1$, $u_2$, $u_3$, $u$.
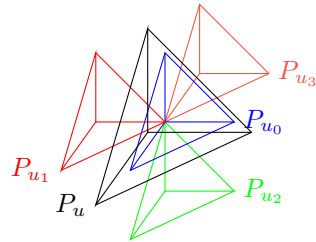
	$$\begin{tikzpicture}[scale=.92]
		\draw [green] (0,0) -- (0,-1) -- (-.5,-1.7) -- (1,-1) -- (0,-1);
		\draw [green] (-.5,-1.7) -- (0,0) -- (1,-1) node [right] {$P_{u_2}$};
		\draw [orange] (0,0) -- (.5,.7) -- (1.5,.7) -- (.5,1.7) -- (.5,.7);
		\draw [orange] (.5,1.7) -- (0,0) -- (1.5,.7)node [right] {$P_{u_3}$};
		\draw [blue] (0,0) -- (1,0)node [right] {$P_{u_0}$} -- (0,1) -- (-.5,-.7) -- (1,0);
		\draw [blue] (0,1)  -- (0,0) -- (-.5,-.7);
		\draw [red] (0,0) -- (-1,0) -- (-1.5,-.7) -- (-1,1) -- (-1,0);
		\draw [red] (-1,1) -- (0,0) -- (-1.5,-.7) node[left] {$P_{u_1}$};
		\draw (1.25,-.15) -- (-.25,-.15) -- (-1,-1.2) -- (-.25,1.35) -- (-.25,-.15);
		\draw (-1,-1.2) node [left] {$P_u$} -- (1.25,-.15) -- (-.25,1.35);
	\end{tikzpicture} $$
	but there is no compatible basis for any maximal cone $\sigma\in\Sigma(3)$.
\end{ex}

\begin{rem}\label{virtual}
	For any equivariant reflexive sheaf $\cE$, we can define a virtual preparliament of polytope. It means that each polytope in the preparliament is considered virtual\footnote{One could also use Cartier divisors instead of virtual polytopes}. This way, we are keeping trace of the defining hyperplanes even if $\cE$ does not have any positivity property. 
\end{rem}

\section{Stability of equivariant vector bundles}\label{s:stability}

Most positivity properties of a toric vector bundle can be detected on its parliaments of polytopes. We may expect other properties such as stability to be visualizable on the parliaments. The problem of classification of vector bundles, or of construction of moduli spaces of semistable coherent sheaves, is still completely open, even in the toric case. 
In this Chapter, we work on finding a necessary and sufficient condition for the (semi)stability of a toric vector bundle $\cE$ (or more generally an equivariant reflexive sheaf) in terms of its virtual (pre)parliament of polytopes $PP_\cE$. This condition can be either seen as an algorithm to check stability or as a visual property on the preparliament $PP_\cE$. In this section we assume, for simplicity, that all the preparliaments are virtual.

\subsection{Definitions}

We start by the definition of slope for coherent sheaves on smooth complete toric varieties. It leads us to the definition of stability for toric vector bundles. The slope of an coherent sheaf $\cE$ depends on the choice of a polarization (an ample divisor $H$ or more generally a movable curve $\alpha$) of the variety $X$ that we leave aside for now and on which we go back in the next subsection.

\begin{defn}
	Let $X$ be a smooth complete toric variety of dimension $d$ and let $\alpha$ be a polarization of $X$. The slope of any coherent sheaf $\cE$ with respect to $\alpha$ is 
	$$\mu_\alpha(\cE):=\frac{c_1(\cE)\cdot\alpha}{\rk(\cE)}\in\bQ\ ,$$
	where $c_1$ is the first Chern class of $\cE$. The dependence on $\alpha$ is often omitted.
\end{defn}

\begin{defn}\label{stabdef}
	Any coherent sheaf $\cE$ is $\alpha$-semistable if for any nonzero subsheaf $\cF\subseteq\cE$, the respective slopes satisfy the inequality $$\mu_\alpha(\cF)\leq\mu_\alpha(\cE)\ .$$ 
	It is $\alpha$-stable if, in addition, for any nonzero subsheaf $\cF\subsetneq\cE$, the strict inequality $\mu_\alpha(\cF)<\mu_\alpha(\cE)$ holds.
	\\
	An equivariant bundle $\cE$ is $\alpha$-polystable if it is a direct sum of $\alpha$-stable bundles.
\end{defn}

\begin{rem}\label{saturated}
	As explained in Remark 2.4 of \cite{HNS}, to prove that some equivariant bundle $\cE$ is $\alpha$-semistable, it is sufficient to verify that $\mu_\alpha(\cF)<\mu_\alpha(\cE)$ holds for every equivariant saturated\footnote{A subsheaf $\cF$ of a torsion-free sheaf $\cE$ is called saturated if $\cE/\cF$ is torsion-free.} subsheaf $\cF$.
\end{rem}

We will need the following combinatorial description of equivariant saturated subsheaves of $\cE$.

\begin{thm}[Proposition 2.3 of \cite{HNS} or Corollary 0.0.2. \cite{DDK}]\label{ers}
	Let $\cE$ be a rank-$r$ equivariant reflexive sheaf on a smooth complete toric variety $X$. Via the Klyachko classification, $\cE$ corresponds to a $\bZ$-filtration $(E^i(j))_{j\in\bZ}$ of $E\cong\bC^r$ for each ray $\rho_i\in\Sigma(1)$.
	
	The equivariant saturated subsheaves $\cF$ of $\cE$ are then in one-to-one correspondence
	with the subfiltrations
	\begin{center}
		$(F^i(j)=E^i(j)\cap F)_{j\in\bZ}\quad$ of $\quad(E^i(j))_{j\in\bZ}$\ ,
	\end{center} 
	for some vector subspace $F\subseteq E$.
\end{thm}

\begin{rem}
	In particular, an equivariant saturated subsheaf $\cF$ of an equivariant reflexive sheaf $\cE$ is reflexive. We thus may construct its preparliament of polytopes as in Definition \ref{ppp}.
\end{rem}

\subsection{Polarization}

Usually, the slope of a coherent sheaf on $X$ is defined with respect to some ample divisor $H$ called a polarization of $X$. More precisely, it involves the self-intersection product $H^{d-1}$. In \cite{GKP}), in order to do birational geometry, Greb, Kebekus and Peternell generalize slopes to movable curves on any complex, projective manifold $X$.

Let us explain it in the toric case. The set of polarizations of $X$ can be extended from ample divisors $H$ to movable divisors $L$ on $X$. This way the movable divisor $L$ polarizes $X$ as well as the blow up of $X$ making $L$ ample. Now, consider a movable curve class $\alpha$ and require it to be positive along a spanning set of rays. It follows from Theorem 3.12, Lemma 4.1 and Theorem 4.2 of \cite{LX} that $\alpha$ can be written as $$\alpha=\langle L_\alpha^{d-1}\rangle $$ where $L_\alpha$ is a unique big and movable divisor and $\langle \rangle $ is the Boucksom positive product\footnote{developed in \cite{BDPP} and constructed algebraically in \cite{BFJ}}. To compute the positive intersection product of big toric divisor classes $L_1$,..., $L_p$, we consider a higher birational model that makes $L_1$,...,$L_p$ nef when removing some positive linear combination of exceptional divisors $D_i$, then $\langle L_1, ..., L_p\rangle$ is defined as 
$$(L_1 - D_1) \cdot ... \cdot (L_p - D_p)\ .$$ 
As such the positive intersection product restricts to the intersection product on ample classes. It justifies that in the following definition of slope, $H^{d-1}$ may be replaced by a movable curve $\alpha$ positive along a spanning set of rays (that is to say $\alpha\cdot D_i>0$ for rays' generators $v_i\in\Sigma(1)$ generating $N_\bR$). 

\begin{defn}\label{polar}
	Let $X$ be a smooth complete toric variety with fan $\Sigma$. 
	A polarization on $X$ is a movable curve $\alpha$ positive along a spanning set of rays of $\Sigma$.
\end{defn}

By Minkowski's theorem for polytopes (see \cite{S} Section 7), the data of a polarization $\alpha$ is equivalent to the data of its weights, defined as follows. 

\begin{defn}\label{weight}
	A movable curve $\alpha$ on $X$, positive along a spanning set of rays, can be written as the positive self-intersection of a unique big and movable divisor $L_\alpha$
	$$\alpha=\langle L_\alpha^{d-1}\rangle\ .$$ Let us denote $P_{L_\alpha}$ the moment polytope of $L_\alpha$ and write $f_i$ for the volume of the face of $P_{L_\alpha}$ which has external normal vector $v_i$. We call the weights of $\alpha$ the numbers $$\text{for }\rho_i\in\Sigma(1),\quad t_i=f_i\frac{(d-1)!}{||v_i||}\geq 0\ .$$ 
\end{defn}

\begin{ex}\label{alpha}
	On $\bP^2$, the polarization $\alpha=D_0=\langle D_0^{2-1}\rangle$ has weights $(1,1,1)$. 
	$$
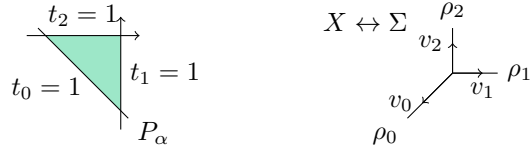
\begin{tikzpicture}[scale=.5]
		\fill [vert!40] (0,0) -- (-2,0) -- (0,-2) -- cycle; 
		\draw [->] (-2.5,0) -- (.5,0);
		\draw [<-] (0,.5) -- (0,-2.5);
		\draw (-2.2,.2) -- (.2,-2.2);
		\draw (.2,-2.6) node [right] {$P_\alpha$};
		\draw (-1,0) node [above] {$t_2=1$};
		\draw (0,-1) node [right] {$t_1=1$};
		\draw (-.8,-1.4) node [left] {$t_0=1$};
	\end{tikzpicture}\qquad \qquad
	\begin{tikzpicture}[scale=.4]
		\draw (-2.9,1.7) node {$X\leftrightarrow\Sigma$};
		\draw (0,0) -- (1.5,0) node[right] {$\rho_1$};
		\draw (0,0) -- (0,1.5)node[above] {$\rho_2$};
		\draw (-1.5,-1.5) node [below left]  {$\rho_0$}-- (0,0);
		\draw [->] (0,0) -- (1,0) node[below] {$v_1$};
		\draw [->] (0,0) -- (0,1) node[left] {$v_2$};
		\draw [<-] (-1,-1) node [left]  {$v_0$} -- (0,0);
	\end{tikzpicture}\vspace{-14pt}$$
\end{ex}

\begin{ex}\label{blowup}
	The polarization $\alpha=2D_1-D_3$ on $\textup{Bl}_{[0:1:0]}\bP^2$ has weights $(1,2,1,1)$.
	
	$$
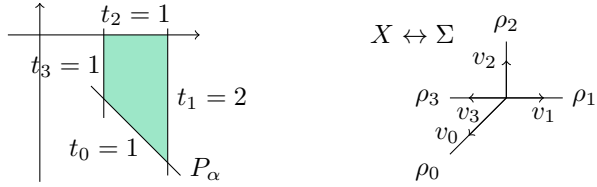
\begin{tikzpicture}[scale=.85]
		\fill [vert!40] (1,0) -- (2,0) -- (2,-2) -- (1,-1); 
		\draw [->] (-.5,0) -- (2.5,0);
		\draw [<-] (0,.5) -- (0,-2.3);
		\draw (.8,-.8) -- (2.2,-2.2);
		\draw (2,.1) -- (2,-2.2);
		\draw (1,.1) -- (1,-1.3);
		\draw (2.2,-2.1) node [right] {$P_\alpha$};
		\draw (1.5,0) node [above] {$t_2=1$};
		\draw (2,-1) node [right] {$t_1=2$};
		\draw (1.1,-.5) node [left] {$t_3=1$};
		\draw (1.7,-1.8) node [left] {$t_0=1$};
	\end{tikzpicture}\qquad \qquad
	\begin{tikzpicture}[scale=.5]
		\draw (-2.5,1.7) node {$X\leftrightarrow\Sigma$};
		\draw (-1.5,0) node[left] {$\rho_3$} -- (1.5,0) node[right] {$\rho_1$};
		\draw (0,0) -- (0,1.5)node[above] {$\rho_2$};
		\draw (-1.5,-1.5) node [below left]  {$\rho_0$}-- (0,0);
		\draw [->] (0,0) -- (1,0) node[below] {$v_1$};
		\draw [->] (0,0) -- (0,1) node[left] {$v_2$};
		\draw [->] (0,0) -- (-1,0) node[below] {$v_3$};
		\draw [<-] (-1,-1) node [left]  {$v_0$} -- (0,0);
	\end{tikzpicture}\vspace{-14pt}$$
\end{ex}

\subsection{Slope and weights}

In this section, we consider a rank-$l$ equivariant saturated subsheaf $\cF$ of an equivariant reflexive sheaf $\cE$ on a toric variety $X$ defined by the filtrations
$$F^i(j)=\left\{\begin{array}{ll}
	F &\text{if }j\leq A^i_{1}\\
	\ldots &\\
	\Span(u_i) &\text{if }A^i_{l-1}<j\leq A^i_{l}\\
	\{0\} &\text{if } A^i_{l}<j.
\end{array}\right.\ .$$ 

We reformulate the slope of $\cF$.
In his original paper \cite{K}, in Remark 3.2.4, Klyachko expressed the Chern classes of any equivariant vector bundle\footnote{See also Proposition 3.1 of \cite{P}}. Kool extended this result, computing the first Chern class of any equivariant coherent sheaf (see Corollary 3.18 of \cite{Ko}). We may reformulate it as follows.

\begin{thm}
	The first Chern class of $\cF$ is given by
	$$c_1(\cF)=\sum_{\rho_i\in\Sigma(1)}\left(\sum_{k=1}^lA_k^i\right) D_i\ .$$
\end{thm}

\begin{prop}\label{c1}
	Let $\cF$ be a rank-$l$ equivariant saturated subsheaf of  $\cE$ as before. If $\alpha$ has weights $(t_i)_{\rho_i\in\Sigma(1)}$ (see Definition \ref{weight}), then the slope of $\cF$ is
	$$\mu_\alpha(\cF)=\sum_i \left(\frac{1}{l}\sum_{k=1}^l A^i_{k}\right) t_i\ .$$
\end{prop}

\begin{prr}
	In Section 4 of \cite{LX}, Lehmann and Xiao show that the numbers $t_i$ are the intersection numbers $D_i\cdot\alpha$.
	We then have that
	$$\vspace{-15pt}\rk(\cF)\mu_\alpha(\cF)=c_1(\cF)\cdot\alpha=\sum_i \left(\sum_{k=1}^l A^i_{k}\right)D_i\cdot\alpha=\sum_i \left(\sum_{k=1}^l A^i_{k}\right) t_i\ .$$
\end{prr}

\subsection{The average polytope}\label{average polytope}

In this section, we consider a nonzero equivariant saturated subsheaf $\cF$ of an equivariant reflexive sheaf $\cE$ on a toric variety $X$. By Theorem \ref{ers}, it is defined by filtrations
$$F^i(j)=\left\{\begin{array}{ll}
	F &\text{if }j\leq A^i_{1}\\
	\ldots &\\
	\Span(u_i)  &\text{if }A^i_{l-1}<j\leq A^i_{l}\\
	\{0\} &\text{if } A^i_{l}<j.
\end{array}\right.\ $$ 
where $F^i(j)=E^i(j)\cap F$, the $(E^i(j))_{\rho_i\in\Sigma(1),j\in\bZ}$ being the defining filtrations of $\cE$.
The slope of $\cF$ may be visualized using the notion of average polytope.

\begin{defn}\label{average}
	We define the average polytope $P_\cF$ of $\cF$ to be the moment polytope associated to the divisor $$\frac{c_1(\cF)}l=\sum_i \left(\frac{1}{l}\sum_{k=1}^l A^i_{k}\right) D_i\ .$$
\end{defn}

\begin{nota}
	We may also use the notation $P_F$ if it appears to be more convenient.
\end{nota}

The average polytope $P_\cF$ of $\cF$ is a visualization of the slopes of $\cF$. In the following proposition, we recover, for any polarization $\alpha$ of $X$, the $\alpha$-slope of $\cF$, only by means of the average polytope $P_\cF$.

\begin{prop}\label{slope,average}
	The $\alpha$-slope of an equivariant saturated subsheaf $\cF\subset\cE$ with average polytope $P_\cF=P_D$ where $D=\sum_{\rho_i\in\Sigma(1)} a_iD_i$ is 
	$$\mu_\alpha(\cF)=\sum_{\rho_i\in\Sigma(1)} a_it_i\ ,$$
	where $\alpha$ has weights $(t_i)_{\rho_i\in\Sigma(1)}$.
	
	Geometrically, the $\alpha$-slope is the mixed volume of the average polytope $P_\cF$ with the polytope $P_{L_\alpha}$ encoding the polarization $$\mu_\alpha(\cF)=V(P_\cF,\underbrace{P_{L_\alpha},...,P_{L_\alpha}}_{d-1\text{ times}})\footnote{This beautiful remark has been pointed out to me by Christian Haase.}\ .$$
\end{prop}

\begin{prr}
	See Proposition \ref{c1} and Definition \ref{average}.
\end{prr}

\begin{ex}
	The average polytopes of $$\cF=\cO(D_0)\oplus\cO(D_2),\quad\cF=\cT_{\bP^2}\quad\text{and}\quad\cF=\cO(D_0)\oplus\cO(D_1)\oplus\cO(D_2)$$ are the hashed polytopes in the following pictures
	$$
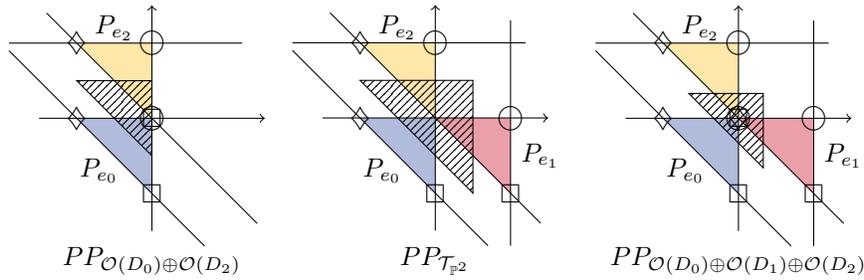
\begin{tikzpicture}[scale=1]
		\fill [jaune!40] (-1,1) -- (0,1) -- (0,0);
		\fill [bleu!40] (-1,0)-- (0,-1) -- (0,0);
		\draw (-1.9,1) -- (1.2,1);
		\draw [->] (-1.5,0) -- (1.5,0);
		\draw [->] (0,-1.5) -- (0,1.5);
		\draw (-1.4,1.4) -- (1.4,-1.4);
		\draw (-1.9,.9) -- (.7,-1.7);    
		\draw (-.5,1.2) node {$P_{e_2}$}; 
		\draw (-.7,-.7) node {$P_{e_0}$};
		\draw (0,-1.9) node {$PP_{\cO(D_0)\oplus\cO(D_2)}$};
		\draw (-1,1) node {$\diamondsuit$};
		\draw (-1,0) node {$\diamondsuit$};
		\draw (0,0) node {$\square$};
		\draw (0,0) node {$\bigcirc$};
		\draw (0,-1) node {$\square$};
		\draw (0,1) node {$\bigcirc$};
		\draw [pattern=north east lines] (0,.5) -- (-1,.5) -- (0,-.5) -- cycle;
	\end{tikzpicture}\quad 
	\begin{tikzpicture}[scale=1]
		\fill [jaune!40] (-1,1) -- (0,1) -- (0,0);
		\fill [bleu!40] (-1,0)-- (0,-1) -- (0,0);
		\fill [rouge!40] (1,-1) -- (1,0) -- (0,0);
		\draw (-1.9,1) -- (1.2,1);
		\draw [->] (-1.5,0) -- (1.5,0);
		\draw (1,-1.7) -- (1,1.3);
		\draw [->] (0,-1.5) -- (0,1.5);
		\draw (-1.4,1.4) -- (1.4,-1.4);
		\draw (-1.9,.9) -- (.7,-1.7);    
		\draw (-.5,1.2) node {$P_{e_2}$};   
		\draw (1.4,-.5) node {$P_{e_1}$};
		\draw (-.7,-.7) node {$P_{e_0}$};
		\draw (0,-1.9) node {$PP_{\cT_{\bP^2}}$};
		\draw (-1,1) node {$\diamondsuit$};
		\draw (-1,0) node {$\diamondsuit$};
		\draw (0,-1) node {$\square$};
		\draw (1,-1) node {$\square$};
		\draw (1,0) node {$\bigcirc$};
		\draw (0,1) node {$\bigcirc$};
		\draw [pattern=north east lines] (.5,.5) -- (-1,.5) -- (.5,-1) -- cycle;
	\end{tikzpicture}\quad
	\begin{tikzpicture}[scale=1]
		\fill [jaune!40] (-1,1) -- (0,1) -- (0,0);
		\fill [bleu!40] (-1,0)-- (0,-1) -- (0,0);
		\fill [rouge!40] (1,-1) -- (1,0) -- (0,0);
		\draw (-1.9,1) -- (1.2,1);
		\draw [->] (-1.5,0) -- (1.5,0);
		\draw (1,-1.7) -- (1,1.3);
		\draw [->] (0,-1.5) -- (0,1.5);
		\draw (-1.4,1.4) -- (1.4,-1.4);
		\draw (-1.9,.9) -- (.7,-1.7);      
		\draw (-.5,1.2) node {$P_{e_2}$};   
		\draw (1.4,-.5) node {$P_{e_1}$};
		\draw (-.7,-.7) node {$P_{e_0}$};
		\draw (0,-1.9) node {$PP_{\cO(D_0)\oplus\cO(D_1)\oplus\cO(D_2)}$};
		\draw (-1,1) node {$\diamondsuit$};
		\draw (-1,0) node {$\diamondsuit$};
		\draw (0,-1) node {$\square$};
		\draw (1,-1) node {$\square$};
		\draw (1,0) node {$\bigcirc$};
		\draw (0,1) node {$\bigcirc$};
		\draw (0,0) node {$\bigcirc$};
		\draw (0,0) node {$\diamondsuit$};		
		\draw (0,0) node {$\square$};
		\draw [pattern=north east lines] (.33,.33) -- (-.66,.33) -- (.33,-.66) -- cycle;
	\end{tikzpicture}$$
	
	The difference between the two last pictures is the following. The hyperplanes forming a '$+$' in the middle have different multiplicities ($1$ or $2$). Indeed, the tangent bundle $\cT_{\bP^2}$ has rank $l=2$ whereas the splitting bundle $\cO(D_0)\oplus\cO(D_1)\oplus\cO(D_2)$ has rank $l=3$. The slopes with respect to the polarization $\alpha$ with weights $(1,1,1)$ (see Example \ref{alpha}) are 
	$$\left\{\begin{array}{l}
		\mu_\alpha\left(\cO(D_0)\oplus\cO(D_2)\right)=\frac12+0+\frac12=1\ ;\\[4pt]
		\mu_\alpha\left(\cT_{\bP^2}\right)=\frac12+\frac12+\frac12=\frac32\ ;\\[4pt]
		\mu_\alpha\left(\cO(D_0)\oplus\cO(D_1)\oplus\cO(D_2)\right)=\frac13+\frac13+\frac13=1\ .
	\end{array}\right.\vspace{-13pt}$$
\end{ex}

\subsection{Result}

In this subsection, we work on a necessary and sufficient condition for stability of equivariant reflexive sheaves. It involves comparing average polytopes by means of the following order. 

\begin{defn}\label{order}
	Let $X$ be a smooth complete toric variety and $\alpha$ a polarization with weights $$t_i:=D_i\cdot \alpha=f_i\frac{(d-1)!}{||v_i||}\quad\text{for every ray $\rho_i\in\Sigma$}\ .$$ Let us denote $P_{L_\alpha}$ the polarization polytope (see Definition \ref{weight}).
	We define a total order on moment polytopes of $X$ by saying
	$$P_1<_\alpha P_2\;\;\iff\;\;\sum_i a^i_{1}t_i<\sum_i a^i_{2}t_i \;\;\left(\text{or } \iff \;\;V(P_1,P_{L_\alpha}^{d-1})<V(P_2,P_{L_\alpha}^{d-1})\right)\ ,$$
	where, for $j\in\{1,2\}$, $P_j$ is the moment polytope of the divisor $\sum_i a^i_{j}D_i$, and $V(\cdot,...,\cdot)$ is the mixed volume of polytopes. 
\end{defn}

Using Definition \ref{order}, Theorem \ref{ers} and Remark \ref{saturated}, we may reformulate the definition of stability (Definition \ref{stabdef}) of toric vector bundles into the following proposition.

\begin{prop}
	An equivariant reflexive sheaf $\cE$ is $\alpha$-stable (resp. $\alpha$-semistable) if and only if for any nonzero subspace $F\subsetneq E$, we have $$P_F<_\alpha P_E \quad(\text{resp. }P_F\leq_\alpha P_E)\ ,$$
	where $P_F$ is the moment polytope of the preparliament $PP_\cF$ obtained by considering the preparliament $PP_\cE$ and keeping only the hyperplanes corresponding to the filtrations \begin{center}$(E^i(j)\cap F)_{i,j}$ of $\cF$\ .\end{center}
\end{prop} 

The following result is a first step towards an algorithm to verify stability.

\begin{thm}\label{thm:existence}
	Let $\cE$ be any equivariant reflexive sheaf on a toric variety $X$. There exists a finite matroid $M(\cE)^S$ representable in $E$ such that $\cE$ is $\alpha$-stable (resp. $\alpha$-semistable) if and only if $$P_F<_\alpha P_E \quad(\text{resp. }P_F\leq_\alpha P_E)\ ,$$
	for any space $\{0\}\subsetneq F=\Span(f)\subsetneq E$ given by a flat $f$ of $M(\cE)^S$.
\end{thm}

\begin{pro}
	Recall that to define (pre)parliaments of polytopes, we introduced, from the filtrations of $\cE$, a semilattice 
	$$L(\cE)=\left\{\bigcap_{i\in\{0,\hdots,n\},}E^i(j_i)\,\middle|\;(j_i)_{i\in\{0,\hdots, n\}}\in\bZ^{n+1}\right\}\ .$$
	Consider an integer $l\in\{1,...,r-1\}$ and let us define the noncontinuous function
	$$\varphi_l: G(l,r)\to \bZ^{L(\cE)}\ ,\quad F\mapsto (\dim(F\cap V))_{V\in L(\cE)}\ .$$
	The $\alpha$-slope of a rank-$l$ saturated subsheaf $\cF$ corresponding to a vector space $F$ depends only on the image of $F$ by $\varphi_l$. 
	As a consequence, to ensure the (semi)stability of $\cE$ we only need to compare the slope of one representative $F\in G(l,r)$ for each image $\varphi_l(F)\in\bZ^{L(\cE)}$ with the slope of $\cE$. Since the image $\im(\varphi_l)\subseteq\{0,...,r\}^{L(\cE)}$ is finite, we have the finite family of subspaces checking stability. We obtain a ground set of a finite matroid computing stability by taking a basis of each member of the family.
\end{pro}

\begin{rem}
	To be more precise, we may define a partial order on $G(l,r)$ by 
	$$F_1<_\cE F_2 \quad\Leftrightarrow \quad \varphi_l(F_1)_V<\varphi_l(F_2)_V\quad\text{for all }V\in L(\cE)\ .$$ The equivariant reflexive sheaf $\cE$ is (semi)stable if and only if for all $l\in\{1,...,r-1\}$, for all maximal element of $\im(\varphi_l)$ with representative $F$, we have that
	$$P_F<_\alpha P_E \quad(\text{resp. }P_F\leq_\alpha P_E)\ .$$
\end{rem}

\begin{rem}
	The matroid $M(\cE)^S$ does not depend on $\alpha$.
\end{rem}

\begin{ex}\label{flatsarenotenough}
	In the next subsections, we prove that $M(\cE)^S$ can be taken to be the matroid $M(\cE)$ of any preparliament of $\cE$, if $\cE$ has rank $1$, $2$ or $3$ or if $\cE=\cT_X$ is the tangent bundle of a smooth complete toric variety. It is not the case in general.
	
	Indeed, we claim that the flats of a matroid $M(\cE)$ are not enough to know that $\cE$ is not (semi)stable.
	Consider $X=\bP^1\times\bP^1$ and the filtrations 
	$$\text{for }\rho_i\in\Sigma(1),\;E^i(j)=\left\{\begin{array}{ll}
		\bC^5 &\text{if }j<0\\
		V_i  &\text{if }0\leq j<1\\
		\{0\} &\text{otherwise}
	\end{array}\right.,$$
	where the vector spaces $V_i$ are defined as
	$V_1=\Span(e_1+e_2,e_3)$, $V_2=\Span(e_1-e_2,e_4)$, $V_3=\Span(e_1,e_5)$ and $V_4=\Span(e_2,e_3+e_4+e_5)$.
	$$
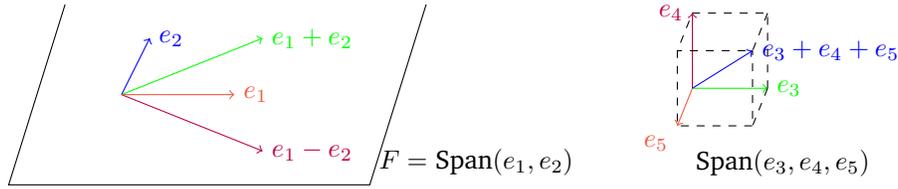
\begin{tikzpicture}[scale=1.5]
		\draw (.5,.8) -- (0,-.8) -- (3.2,-.8) node[above right] {$F=\textup{Span}(e_1,e_2)$} -- (3.7,.8);
		\draw [->,purple] (1,0) -- (2.25,-.5) node[right] {$e_1-e_2$};
		\draw [->,green] (1,0) -- (2.25,.5) node[right] {$e_1+e_2$};
		\draw [->,blue] (1,0) -- (1.25,.5) node[right] {$e_2$};
		\draw [->,orange] (1,0) -- (2,0) node[right] {$e_1$};
	\end{tikzpicture}\qquad
	\begin{tikzpicture}
		\draw [->,green] (0,0) -- (1,0) node [right] {$e_3$};
		\draw [->,purple] (0,0) -- (0,1) node [left] {$e_4$};
		\draw [->,orange] (0,0) -- (-.2,-.5) node [below left] {$e_5$};
		\draw [->,blue] (0,0) -- (.8,.5) node [right] {$e_3+e_4+e_5$};
		\draw [dashed] (1,0) -- (.8,-.5) -- (.8,.5); 
		\draw [dashed] (1,0) -- (1,1) -- (0,1) -- (-.2,.5) -- (-.2,-.5) -- (.8,-.5);
		\draw [dashed] (1,1) -- (.8,.5) -- (-.2,.5);
		\draw (1.2,-1) node {$\textup{Span}(e_3,e_4,e_5)$};
	\end{tikzpicture}$$
	A matroid that can be obtained by Algorithm 3.2	is given by the ground set
	\begin{align*}
		G=\{ &e_1+e_2+e_3,\ e_1+e_2-e_3,\ e_1-e_2+2e_4,\ e_1-e_2-2e_4,\\ &e_1+3e_5,\ e_1-3e_5,\ e_2+4e_3+4e_4+4e_5,\ e_2-4e_3-4e_4-4e_5 \} \ .
	\end{align*}
	Consider the vector space $F=\Span(e_1,e_2)$. The slope of $\cE$ and $\cF$ are 
	$$\mu(\cE)=\frac{2+2+2+2}5=\frac85\quad\text{and}\quad\mu(\cF)=\frac{1+1+1+1}2=2\ .$$
	Consequently $\cE$ is unstable, but it cannot be seen by looking at the flats of $G$.
	Indeed, we have that 
	\begin{itemize}
		\item A flat $f$ of dimension $l=1$, $2$ or $3$ satisfies $$\textup{Span}(V_i\cap f)=V_i\cap\textup{Span}(f)$$ and its associated saturated subsheaf $\cF$ has thus slope $\mu(\cF)=l/l=1$.
		\item A flat of dimension $4$ either $$\textup{Span}(V_i\cap f)=V_i\cap\textup{Span}(f)$$ and its associated saturated subsheaf has slope $4/4=1$,
		or is of the form $V_i+V_j=\Span(e_1,e_2,e'_1,e'_2)$ with $e'_1,e'_2$ in $G$ and its associated saturated subsheaf $\cF$ has slope $\mu(\cF)=\frac64=\frac32$.
	\end{itemize} \vspace{-13pt}
\end{ex}

\begin{rem}
	We may enhance Example \ref{flatsarenotenough} by proving that all the flats of all the matroids $M(\cE)$ of a preparliament are not enough to check the (semi)stability of $\cE$.
	Indeed, we may force the ground set to be exactly the set $G$ of Example \ref{flatsarenotenough}. It suffices to incorporate $1$-dimensional spaces $V=\Span(g)$ for any $g\in G$ into the filtrations (take for instance $X=\bP^1\times\bP^1\times\bP^1\times\bP^1$).
\end{rem}

\subsection{Case of rank less or equal to 3}

The purpose of the following theorem is to prove that if $\cE$ has rank less or equal to $3$, then the matroid $M(\cE)$ of any preparliament of $\cE$ checks stability.

\begin{thm}\label{result}
	Let $X$ be a smooth complete toric variety, let $\cE$ be an equivariant reflexive sheaf on $X$, of rank less or equal to $3$, and let $(PP_\cE,G(\cE))$ be any of its preparliament. Then $\cE$ is $\alpha$-(semi)stable if and only if for any nonzero flat $f\subsetneq G(\cE)$, we have $$P_\cF<_\alpha P_{G(\cE)}\quad(\text{resp. }P_\cF\leq_\alpha P_{G(\cE)})\ ,$$
	where $\cF$ is the equivariant saturated sheaf corresponding by Theorem \ref{ers} to the linear subspace $\Span(f)\subset E$.
\end{thm}

Before proving the theorem, let us start by a lemma.

\begin{lem}\label{slope}
	Consider two vector spaces $F$ and $F_0$ of the same dimension $l$. Assume that for every ray $\rho_i\in\Sigma(1)$, for any $j\in\bZ$, we have that
	$$\dim(F\cap E^i(j))\geq\dim(F_0\cap E^i(j))\ .$$ Then the slope of $\cF$ is greater than or equal to the slope of $\cF_0$
	$$\mu(\cF)\geq\mu(\cF_0)\ .$$
\end{lem}

\begin{prr}
	If $A^i_k$ and $B^i_k$ are the remaining gaps after intersecting the filtration $(E^i(j))_{j\in\bZ}$ respectively with $F_0$ and $F$ then for all $k\in\{1,\hdots,l\}$, $\rho_i\in\Sigma(1)$,
	$$\left\{\begin{aligned}
		A^i_k&=\min\left\{j\ |\ \dim(F_0\cap E^i(j))\leq l-k+1\right\}\\ 
		B^i_k&=\min\left\{j\ |\ \dim(F\cap E^i(j))\leq l-k+1\right\}
	\end{aligned}\right.$$
	satisfy $A^i_k\geq B^i_k$.
	We conclude by Proposition \ref{c1}.
\end{prr}

We may now start the proof of Theorem \ref{result}.

\begin{prr}
	By Lemma \ref{slope}, it is enough to show that for any vector space $\{0\}\subsetneq F_0\subsetneq E$, there exists a flat $f\subsetneq G(\cE)$ with $F=\Span(f)$ of the same dimension than $F_0$ and satisfying 
	$$\dim(F\cap V)\geq\dim(F_0\cap V)\quad\text{for all }V\in L(\cE)\ .$$
	
	In the case where $\dim(F_0)=1$, take $F$ generated by a vector of $G(\cE)$ in $\bigcap_{V\text{ st }F_0\cap V\neq\{0\}} V$. It is possible because $\bigcap_{F_0\cap V\neq\{0\}} V$ is nonzero ($F_0$ being included into it) and $G(\cE)$ has been constructed such that $\Span(G(\cE)\cap \bigcap_{V\text{ st }F_0\cap V\neq\{0\}} V)=\bigcap_{V\text{ st }F_0\cap V\neq\{0\}} V$. It concludes the proof for bundles of rank $2$.
	
	It remains the case where $\dim(F_0)=2$ and $\rk(\cE)=3$. Let us denote by $$V_1^1,...,V_a^1,V_1^2,...,V_b^2$$ the vector spaces of $L(\cE)$ such that $\dim(F_0\cap V_k^\varepsilon)=\varepsilon$. Remark that $b$ is either equal to $1$ (and $V_1^2=E$) or $2$ (and $\{V_1^2,V_2^2\}=\{F,E\}$). If $b=2$ then there exists a flat $f\subseteq G(\cE)$ such that $\Span(f)=F_0$ and we are done. Assume thus that $b=1$ and let us distinguish cases.
	\begin{itemize}
		\item If there exists $k_1,k_2\in\{1,...,a\}$ such that $V_{k_1}^1,V_{k_2}^1$ of dimension $1$, then $F_0$ is also generated by a flat $$f=\{e_1,e_2\}$$ of $G(\cE)$, where $V_{k_1}^1=\Span(e_1),\ V_{k_2}^1=\Span(e_2)$. We can take $F=F_0$.
		\item If there exists exactly one $V_{k_1}^1$ of dimension $1$, then take a generator $e_1\in G(\cE)$ of $V_k^1$ and any other vector $e_2\in G(\cE)$. We define $F$ as $F=\Span(e_1,e_2)$. As for all $k\neq k_1$, we have that $\dim(V_{k}^1)\geq2$, $\dim(F)=2$ and both $V_{k}^1$ and $F$ are included in $E$ of dimension $3$, it holds that $$\dim(F\cap V_k^1)\geq1=\dim(F_0\cap V_k^1)\ .$$
		\item If for all $k\in\{1,...,a\}$, we have $\dim(V_k^1)\geq2$ then, by our last argument, any dimension $2$-space $F=\Span(f)$ generated by a flat of $G(\cE)$ satisfies $$\dim(F\cap V_k^1)\geq1=\dim(F_0\cap V_k^1)\quad\text{ for all }k\in\{1,...,a\}\ .$$
	\end{itemize}\vspace{-25pt}
\end{prr}

\subsection{A generalisation of the case of tangent bundles}

Let $X$ be a smooth complete toric variety of dimension $d$ and let $\cE$ be any equivariant reflexive sheaf $\cE$ on $X$ with Klyachko filtrations given by vector spaces of dimension $0$, $1$ and $d$ and general gaps $A^i_j$. For instance, the tangent bundle $\cE=\cT_X$ on $X$ satisfies these condition as its Klyachko description is
$$\text{for }\rho_i\in\Sigma(1),\;E^i(j)=\left\{\begin{array}{ll}
	\bC^d &\text{if }j<0\\
	\Span(v_i)  &\text{if }0\leq j<1\\
	\{0\} &\text{otherwise}
\end{array}\right..$$ 
The tangent bundle case has already been treated in \cite{HNS} or \cite{DDK1}. In this section, we prove that $M(\cE)$ is a matroid that check stability.
The only matroid $M(\cE)$ associated to $\cE$ has ground set given by 
$G(\cE)=\{e_i\ |\ \rho_i\in\Sigma(1)\}$, where the spaces $\Span(e_i)$ are the $1$-dimensional vector spaces appearing in the Klyachko filtrations of $\cE$.

\begin{prop}
	The equivariant reflexive sheaf $\cE$ is $\alpha$-(semi)stable if and only if for any nonzero flat $f\subsetneq G(\cE)$, we have $$P_\cF<_\alpha P_{G(\cE)}\quad(\text{resp. }P_\cF\leq_\alpha P_{G(\cE)})\ ,$$
	where $\cF$ is the equivariant saturated sheaf corresponding by Theorem \ref{ers} to the linear subspace $\Span(f)\subset E$.
\end{prop}

\begin{prr}
	The vector spaces in 	$$L(\cE)=\left\{\bigcap_{i\in\{0,\hdots,n\},}E^i(j_i)\,\middle|\;(j_i)_{i\in\{0,\hdots, n\}}\in\bZ^{n+1}\right\}$$
	are of dimension either $0$, $1$ or $d$.
	Hence, any saturated subsheaf $\cF_0\subset\cE$ corresponds to a subspace $F_0\subset \bC^d$ having intersection with elements of $L(\cE)$ either
	$$F_0\cap V=\{0\}\ ,\qquad F_0\cap V=V \text{ and }\dim(V)=1\ ,\qquad\text{or } F_0\cap V=F_0\ .$$
	Each $V\in L(\cE)$ of dimension $1$ has a generator $g_V$ in $G(\cE)$. Let us consider a flat $f$ of $G(\cE)$ obtained by considering the generators $g_V\in G(\cE)$ of the vector spaces $V\in L(\cE)$ of dimension $1$ included in $F_0$, and by adjoining elements of $G$ such that the vector space $F:=\Span(f)$ satisfies	$\dim(F)=\dim(F_0).$ 
	By Lemma \ref{slope}, we have the inequality
	$$\mu(\cF)\geq\mu(\cF_0)\ ,$$
	and the flats of $G(\cE)$ are enough to check the (semi)stability of $\cE$. 
\end{prr}

\begin{rem}
	In the tangent bundle case, we recover the criterion for stability given in Proposition 1.2. of \cite{HNS}.
\end{rem}

\begin{ex}
	Consider the tangent bundle $\cE=\cT_{\bP^2}$ on $\bP^2$. The nontrivial flats of $G(\cE)=\{v_0,v_1,v_2\}$ are \begin{center}$f_0=\{v_0\}$, $f_1=\{v_1\}$ and $f_2=\{v_2\}$.\end{center} They correspond to the three nontrivial equivariant saturated subsheaves $\cF_0$, $\cF_1$ and $\cF_2$ of $\cE$. Let us look at the average polytopes: $$
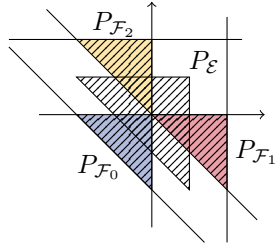
\begin{tikzpicture}[scale=1]
		\fill [jaune!40] (-1,1) -- (0,1) -- (0,0);
		\fill [bleu!40] (-1,0)-- (0,-1) -- (0,0);
		\fill [rouge!40] (1,-1) -- (1,0) -- (0,0);
		\draw (-1.9,1) -- (1.2,1);
		\draw [->] (-1.5,0) -- (1.5,0);
		\draw (1,-1.7) -- (1,1.3);
		\draw [->] (0,-1.5) -- (0,1.5);
		\draw (-1.4,1.4) -- (1.4,-1.4);
		\draw (-1.9,.9) -- (.7,-1.7);    
		\draw [pattern=north east lines] (-1,1) -- (0,1) -- (0,0);
		\draw [pattern=north east lines] (-1,0)-- (0,-1) -- (0,0);
		\draw [pattern=north east lines] (1,-1) -- (1,0) -- (0,0);
		\draw (-.5,1.2) node {$P_{\cF_2}$};   
		\draw (1.4,-.5) node {$P_{\cF_1}$};
		\draw (-.7,-.7) node {$P_{\cF_0}$};
		\draw (.7,.7) node {$P_{\cE}$};
		\draw [pattern=north east lines] (.5,.5) -- (-1,.5) -- (.5,-1) -- cycle;
	\end{tikzpicture}$$	
	By Theorem \ref{result}, the fact that $P_{\cF_0},P_{\cF_1},P_{\cF_2}<_\alpha P_\cE$ proves that $\cE$ is $\alpha$-stable.
\end{ex}

\section{Parliaments of equivariant subbundles}\label{s:sub}

In this section, we first define subparliaments of polytopes and identify them with parliaments of equivariant subbundles. Second, we give an example of toric subbundle which is not a direct factor.

\subsection{Subparliaments as parliaments of toric subbundles}

Let us state our definition of subparliament of a parliament of polytopes.

\begin{defn}\label{com}
	Let $f=G(\cE)\cap F$ be a flat of the matroid $M(\cE)$ such that there exists a compatible basis $(B_\sigma)_{\sigma\in\Sigma(d)}$ for $\cE$ satisfying
	$$\text{for all }\sigma\in\Sigma(d),\quad\#B_\sigma\cap f=\dim F\ .$$
	Then we call $f$ a compatible flat, and we call the subset of polytopes of $PP_\cE$ indexed by elements of $f=F\cap G(\cE)$
	$$\left\{(P_e,e)\ \middle|\ e\in f\right\}$$
	a subparliament of the parliament of polytopes $\cE$.
\end{defn}

Equivariant subbundles have been combinatorially described in \cite{KD}. We translate their results in terms of parliaments of polytopes.

\begin{prop}\label{subb}
	Via the Klyachko classification, a rank-$r$ equivariant vector bundle $\cE$ on $X$ corresponds to a $\bZ$-filtration $(E^i(j))_{j\in\bZ}$ of $E\cong\bC^r$ for each $\rho_i\in\Sigma(1)$.
	
	The equivariant subbundles $\cF$ of $\cE$ are then in one-to-one correspondence with the subfiltrations \begin{center}
		$(F^i(j)=E^i(j)\cap F)_{j\in\bZ}\quad$of$\quad(E^i(j))_{j\in\bZ}$\ ,
	\end{center} 
	for some vector subspace $F\subseteq E$ such that there exist ground sets $G(\cF)$ and $G(\cE)$ resulting from Algorithm 3.2 for $L(\cF)$ and $L(\cE)$ with
	$$G(\cE)\cap F=G(\cF)$$
	and such that there exist compatible bases $(B_\sigma)_{\sigma\in\Sigma(d)}$ for $\cE$ which, restricted to $F$, satisfy the compatibility condition (\ref{CC}) for $\cF$.
\end{prop}

\begin{prr}
	By Proposition 4.1.1 of \cite{KD}, the equivariant subbundles $\cF$ of $\cE$ are then in one-to-one correspondence with the subfiltrations \begin{center}
		$(F^i(j)=E^i(j)\cap F)_{j\in\bZ}\quad$of$\quad(E^i(j))_{j\in\bZ}$\ ,
	\end{center} 
	for some vector subspace $F\subseteq E$ such that $\{F,\{E^i(j)\}_{\rho_i\in\sigma(1)}\}$ of $E$ forms a distributive lattice. 
	
	Let us show that there exist ground sets $G(\cF)$, $G(\cE)$ obtained by Algorithm 3.2 for $L(\cF),L(\cG)$ with
	$$G(\cE)\cap F=G(\cF)\ .$$
	Consider the ground set $G(\cE)$ obtained by taking a maximal number of vectors in $F$ in Algorithm 3.2. 
	If $V_1,V_2\in L(\cE)$, let us say that $V_1<V_2$ if $V_1$ appears before $V_2$ in Algorithm 3.2. We denote by $G_V$ the intersection of $V$ with the temporary ground set at the beginning of Step $V\in L(\cE)$.
	During Step $V$, we need to add a basis $CB_V$ of any complement of $\Span(G_V)$ in $V$ to the ground set. We choose this basis to be composed of a maximal number of vectors in $F$ so that we have
	\begin{align*}
		V\cap F&=\left(\Span(G_V)+\Span(CB_V)\right)\cap F\\
		&=\left(\sum_{\substack{V'=V_1\cap V\\V_1<V}}V'+\Span(CB_V)\right)\cap F\\
		&=\left(\sum_{\substack{V'=V_1\cap V\\V_1<V}}V'\right)\cap F+\Span(CB_V\cap F)\\
		&=\sum_{\substack{V'=V_1\cap V\\V_1<V}}(V'\cap F)+\Span(CB_V\cap F)\ .
	\end{align*}
	The last step comes from $\{F,\{E^i(j)\}_{\rho_i\in\sigma(1)}\}$ forming a distributive lattice of $E$.
	By induction on the $V\in L(\cE)$ appearing in Algorithm 3.2, $G_V\cap F$ generates every $(V_1\cap V)\cap F$ for $V_1<V$. We finally obtain that $G(\cE)\cap F$ generates every $V\cap F$ for $V\in L(\cE)$.
	\\
	Now if we denote by $G_{V\cap F}^\cF$ the intersection of $V\cap F$ with the temporary ground set at the beginning of Step $V\cap F\in L(\cF)$ of Algorithm 3.2 (applied to $\cF$), then
	$$V\cap F=\Span(G_{V\cap F}^\cF)+\Span(CB_{V\cap F}^\cF)=
	\sum_{\substack{V'\cap F=(V_1\cap F)\cap (V\cap F)\\(V_1\cap F)<(V\cap F)}}(V'\cap F)+\Span(CB_{V\cap F}^\cF)\ .$$
	Moreover, since $V_1<V\Rightarrow(V_1\cap F<V\cap F)$ and $(V_1\cap F<V\cap F)\Rightarrow (V_1\cap V<V)$, we have 
	\begin{align*}
		&\left\{V'\cap F\ \middle|\ V'=V_1\cap V\text{ and } V_1<V\right\}\\=\ &\left\{V'\cap F\ \middle|\ V'\cap F=(V_1\cap F)\cap (V\cap F)\text{ and }(V_1\cap F)<(V\cap F)\right\}\ . 
	\end{align*}
	We may take $CB_{V\cap F}^\cF$ to be $CB_V\cap F$, and we finally obtain $G(\cF)=G(\cE)\cap F$.
	
	Moreover, for each cone $\sigma\in\Sigma(d)$, any compatible basis $B_\sigma^F$ of $F$ extends to a compatible basis $B_\sigma$ of $E$. Indeed, to construct a compatible basis $B_\sigma$ (resp. $B_\sigma^F$), we apply Algorithm 3.2 to 
	$$L_\sigma(\cE)=\left\{\bigcap_{\rho_i\subset\sigma}E^i(j_i)\middle|(j_i)_{\rho_i\subset\sigma}\right\}\quad\text{(resp. $L_\sigma(\cF)=\left\{V\cap  F\middle| V\in L_\sigma(\cE)\right\}$)}\ ,$$
	taking $CB_V$ included in $G(\cE)$ (resp. taking $CB_{V\cap F}^\cF$ included in $G(\cE)\cap F$). Extending the bases $B_\sigma^F$ is possible because $G(\cE)$ generates each $V\in L_\sigma(\cE)\subset L(\cE)$ as well as $G(\cF)=G(\cE)\cap F$ generates each $V\cap F\in L_\sigma(\cF)\subset L(\cF)$. 
\end{prr}

\begin{defn}
	Consider a toric vector bundle $\cE$ and a parliament of polytopes $PP_\cE$ of $\cE$. 
	A subparliament of polytopes of $PP_\cE$ is a set of the form
	$$PP_\cF=\left\{(P_e,e)\ \middle|\ e\in G(\cE)\cap F\right\} \subset PP_\cE\ ,$$
	for some vector space $F$, which is a parliament meaning that it satisfies the compatibility condition (\ref{CC}).
\end{defn}

\begin{thm}\label{parl}
	Let $\cE$ be a toric vector bundle and let $\cF$ be the equivariant subbundle of $\cE$ corresponding by Klyachko's theorem to the filtrations $(E^i(j)\cap F)_{j\in\bZ}$ for each $\rho_i\in\Sigma(1)$.
	
	Then, there exists a parliament $PP_\cE$ of $\cE$ such that a parliament of $\cF$ is a subparliament of $P_\cE$
	$$PP_\cF=\left\{(P_e,e)\ \middle|\ e\in G(\cE)\cap F\right\} \subset PP_\cE\ .$$
\end{thm}

\begin{prr}	
	Consider a ground set $G(\cE)$ for the parliament of $\cE$ obtained via Algorithm 3.2 by taking at each step a maximal number of elements in $F$. Then by Proposition \ref{subb}, $G(\cF)=G(\cE)\cap F$ is a possible output for Algorithm 3.2 applied to $L(\cF)$. Now if $e\in G(\cF)=G(\cE)\cap F$ then the polytope $P_e\in PP_\cF$ and
	\begin{align*}
		P_e&=\left\{m\,\middle|\;\forall i,\:\langle m,v_i\rangle \leq\max\left\{j\,\middle|\;e\in E^i(j)\cap F\right\}\right\}\\
		&=\left\{m\,\middle|\;\forall i,\:\langle m,v_i\rangle \leq\max\left\{j\,\middle|\;e\in E^i(j)\right\}\right\}
	\end{align*} 
	is the same as the one in the parliament of $\cE$.
\end{prr}

\begin{ex}
	Consider the toric variety $X=\bP^2$ with its $T$-invariant divisors $D_0$, $D_1$ and $D_2$ and the splitting, equivariant, rank-$2$ vector bundle $\cE=\cO(D_0)\oplus\cO(D_1+D_2)$ on $X$. The corresponding $\bZ$-filtrations of $E=\bC^2$ are of the form $$E^i(j)=\left\{\begin{array}{ll}
		\bC^2 &\text{if }j\leq0\\
		\Span(e_1)  &\text{if }0<j\leq1\\
		\{0\} &\text{otherwise}
	\end{array}\right.\;(\text{for }i=1,2),\quad E^0(j)=\left\{\begin{array}{ll}
		\bC^2 &\text{if }j\leq0\\
		\Span(e_2)  &\text{if }0<j\leq1\\
		\{0\} &\text{otherwise}
	\end{array}\right.$$
	where $\Span(e_1) $ and $\Span(e_2) $ are different lines in $\bC^2$. The ground set is $G(\cE)=\{e_1,e_2\}$. The toric subbundle $\cF=\cO(D_0)\subset\cE$ corresponds to the flat $f=\{e_2\}\subset G(\cE)$. 
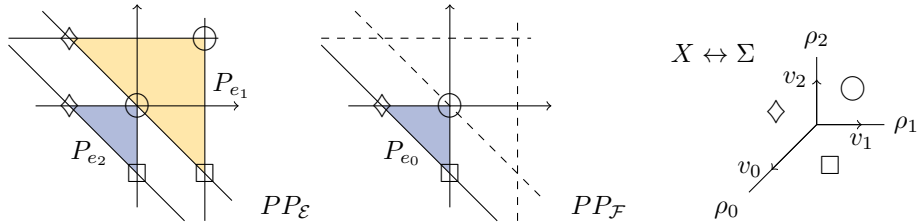
	$$
	\begin{tikzpicture}[scale=.9]
		\fill [jaune!40] (-1,1) -- (1,1) -- (1,-1);
		\fill [bleu!40] (-1,0)-- (0,-1) -- (0,0);
		\draw (-1.9,1) -- (1.2,1);
		\draw [->] (-1.5,0) -- (1.5,0);
		\draw (1,-1.7) -- (1,1.3);
		\draw [->] (0,-1.5) -- (0,1.5);
		\draw (-1.4,1.4) -- (1.4,-1.4);
		\draw (-1.9,.9) -- (.7,-1.7); 
		\draw (-1,1) node {$\diamondsuit$};
		\draw (1,1) node {$\bigcirc$};
		\draw (1,-1) node {$\square$};
		\draw (-1,0)node {$\diamondsuit$};
		\draw (0,-1) node {$\square$};
		\draw (0,0) node {$\bigcirc$};
		\draw (1.4,.3) node {$P_{e_1}$};
		\draw (-.7,-.7) node {$P_{e_2}$};
		\draw (2.2,-1.5) node {$PP_\cE$};
	\end{tikzpicture}
	\begin{tikzpicture}[scale=.9]
		\fill [bleu!40] (-1,0)-- (0,-1) -- (0,0);
		\draw [dashed] (-1.9,1) -- (1.2,1);
		\draw [->] (-1.5,0) -- (1.5,0);
		\draw [dashed] (1,-1.7) -- (1,1.3);
		\draw [->] (0,-1.5) -- (0,1.5);
		\draw [dashed] (-1.4,1.4) -- (1.4,-1.4);
		\draw (-1.9,.9) -- (.7,-1.7);  
		\draw (-.7,-.7) node {$P_{e_0}$};
		\draw (0,-1) node {$\square$};
		\draw (0,0) node {$\bigcirc$};
		\draw (-1,0)node {$\diamondsuit$};
		\draw (2.2,-1.5) node {$PP_\cF$};
	\end{tikzpicture} \quad\begin{tikzpicture}[scale=.6]
		\draw (-2.3,1.6) node {$X\leftrightarrow\Sigma$};
		\draw (0,0) -- (1.5,0) node[right] {$\rho_1$};
		\draw (0,0) -- (0,1.5)node[above] {$\rho_2$};
		\draw (-1.5,-1.5) node [below left]  {$\rho_0$}-- (0,0);
		\draw [->] (0,0) -- (1,0) node[below] {$v_1$};
		\draw [->] (0,0) -- (0,1) node[left] {$v_2$};
		\draw [<-] (-1,-1) node [left]  {$v_0$} -- (0,0);
		\draw (-.9,.3) node {$\diamondsuit$};
		\draw (.3,-.9) node {$\square$};
		\draw (.8,.8) node {$\bigcirc$};
		\draw [white] (0,-2)--(0,-2.2);
	\end{tikzpicture}$$ 
	The compatibility condition (\ref{CC}) for $\cE$ is satisfied on taking $$B_\sigma=\left\{L_u^\sigma\ \middle|\ u\in\textbf u(\sigma)\right\}=\{e_1,e_2\}\quad\text{for all }\sigma\in\Sigma(d) .$$
	We see that for all $\sigma\in\Sigma(d)$, $B_\sigma\cap f=\{e_2\}$ generates $F=\Span(f)$, the flat $f=\{e_2\}$ is thus compatible. The parliament $PP_\cF$ is a subparliament of $PP_\cE$. 
\end{ex}

\begin{ex}\label{tp22}
	Let us come back to Example \ref{tp2} with $X=\bP^2$ and its tangent bundle $\cE=\cT_{\bP^2}$.
	$$
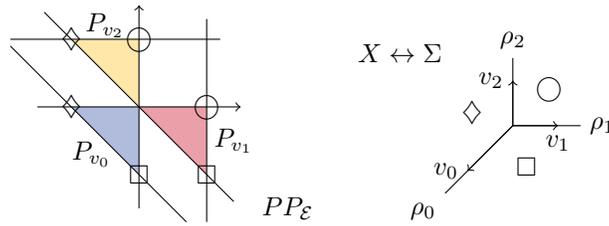
\begin{tikzpicture}[scale=.9]
		\fill [jaune!40] (-1,1) -- (0,1) -- (0,0);
		\fill [bleu!40] (-1,0)-- (0,-1) -- (0,0);
		\fill [rouge!40] (1,-1) -- (1,0) -- (0,0);
		\draw (-1.9,1) -- (1.2,1);
		\draw [->] (-1.5,0) -- (1.5,0);
		\draw (1,-1.7) -- (1,1.3);
		\draw [->] (0,-1.5) -- (0,1.5);
		\draw (-1.4,1.4) -- (1.4,-1.4);
		\draw (-1.9,.9) -- (.7,-1.7);   
		\draw (-.5,1.2) node {$P_{v_2}$};   
		\draw (1.4,-.5) node {$P_{v_1}$};
		\draw (-.7,-.7) node {$P_{v_0}$};
		\draw (2.2,-1.5) node {$PP_\cE$};
		\draw (-1,1) node {$\diamondsuit$};
		\draw (-1,0) node {$\diamondsuit$};
		\draw (0,-1) node {$\square$};
		\draw (1,-1) node {$\square$};
		\draw (1,0) node {$\bigcirc$};
		\draw (0,1) node {$\bigcirc$};
	\end{tikzpicture}\quad \begin{tikzpicture}[scale=.6]
		\draw (-2.5,1.6) node {$X\leftrightarrow\Sigma$};
		\draw (0,0) -- (1.5,0) node[right] {$\rho_1$};
		\draw (0,0) -- (0,1.5)node[above] {$\rho_2$};
		\draw (-1.5,-1.5) node [below left]  {$\rho_0$}-- (0,0);
		\draw [->] (0,0) -- (1,0) node[below] {$v_1$};
		\draw [->] (0,0) -- (0,1) node[left] {$v_2$};
		\draw [<-] (-1,-1) node [left]  {$v_0$} -- (0,0);
		\draw (-.9,.3) node {$\diamondsuit$};
		\draw (.3,-.9) node {$\square$};
		\draw (.8,.8) node {$\bigcirc$};
	\end{tikzpicture}$$
	There exists no nontrivial equivariant subbundle $\cF$ of $\cE$ because the nontrivial flats are $f_0=\{v_0\}$, $f_1=\{v_1\}$ and $f_2=\{v_2\}$ and are not compatible.
\end{ex}

\begin{cor}\label{subp}
	The subparliaments of the parliaments of $\cE$ are the parliaments of the equivariant subbundles of $\cE$.
\end{cor}

\subsection{Subbundle that is not a direct factor}

There exists equivariant vector bundles $\cE$ with a nontrivial equivariant subbundle $\cF$ that we cannot factorize in a direct sum of $\cE$. 

\begin{ex}\label{nontriv}
	Consider the equivariant vector bundle $\cE$ on $\bP^2$ defined by the filtrations 
	\begin{align*}&E^0(j)=\left\{\begin{array}{ll}
			\bC^3 &\text{if }j<-3\\
			\Span(e_1,e_2) &\text{if }-3\leq j<-1\\
			\Span(e_1)  &\text{if }-1\leq j<1\\
			\{0\} &\text{otherwise}
		\end{array}\right.,\quad
		E^1(j)=\left\{\begin{array}{ll}
			\bC^3 &\text{if }j<0\\
			\Span(e_2,e_3) &\text{if }0\leq j<2\\
			\Span(e_3)  &\text{if }2\leq j<4\\
			\{0\} &\text{otherwise}
		\end{array}\right.\\
		&\text{and}\quad
		E^2(j)=\left\{\begin{array}{ll}
			\bC^3 &\text{if }j<0\\
			\Span(e_1-e_3,e_1-e_2) &\text{if }0\leq 	j<2\\
			\Span(e_1-e_3) &\text{if }2\leq j<4\\
			\{0\} &\text{otherwise}
		\end{array}\right.\end{align*} 
	where $e_1=(1,0,0)$, $e_2=(0,1,0)$, $e_3=(0,0,1)$.
	The unique possible ground set of $M(\cE)$ is $$G(\cE)=\{e_1,e_2,e_3,e_1-e_2,e_1-e_3,e_2-e_3\}\ .$$
	Taking $F=\Span(e_2)^\perp$, we obtain an equivariant subbundle $\cF$. 
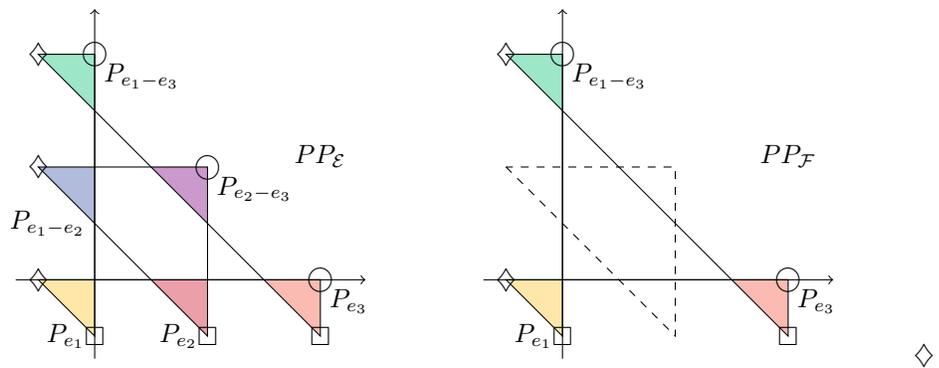
	$$
	\begin{tikzpicture}[scale=1.5]
		\fill [jaune!40] (-.5,0) -- (0,-.5) -- (0,0);
		\fill [bleu!40] (-.5,1) -- (0,.5) -- (0,1);
		\fill [vert!40] (-.5,2) -- (0,1.5) -- (0,2);
		\fill [rouge!40] (.5,0) -- (1,-.5) -- (1,0);
		\fill [orange!40] (1.5,0) -- (2,-.5) -- (2,0);
		\fill [violet!40] (.5,1) -- (1,.5) -- (1,1);
		\draw [->] (-.7,0) -- (2.4,0);
		\draw [->] (0,-.7) -- (0,2.4);
		\draw (-.5,0) -- (2,0) -- (2,-.5) -- (-.5,2) -- (0,2) -- (0,-.5) -- cycle; 
		\draw (-.5,1) -- (1,-.5) -- (1,1) -- cycle;
		\draw (0,-.5) node [left] {$P_{e_1}$};
		\draw (0,.5) node [left] {$P_{e_1-e_2}$};
		\draw (0,2) node [below right] {$P_{e_1-e_3}$};
		\draw (1,-.5) node [left] {$P_{e_2}$};
		\draw (2,0) node [below right] {$P_{e_3}$};
		\draw (1,1) node [below right] {$P_{e_2-e_3}$};
		\draw (-.5,0) node {$\diamondsuit$};
		\draw (-.5,1) node {$\diamondsuit$};
		\draw (-.5,2) node {$\diamondsuit$};
		\draw (0,-.5) node {$\square$};
		\draw (1,-.5) node {$\square$};
		\draw (2,-.5) node {$\square$};
		\draw (2,0) node {$\bigcirc$};
		\draw (1,1) node {$\bigcirc$};
		\draw (0,2) node {$\bigcirc$};
		\draw (2,1.1) node {$PP_{\cE}$};
	\end{tikzpicture}\qquad\qquad\vspace{-15pt}
	\begin{tikzpicture}[scale=1.5]
		\fill [vert!40] (-.5,2) -- (0,1.5) -- (0,2);
		\fill [jaune!40] (-.5,0) -- (0,-.5) -- (0,0);
		\fill [orange!40] (1.5,0) -- (2,-.5) -- (2,0);
		\draw [->] (-.7,0) -- (2.4,0);
		\draw [->] (0,-.7) -- (0,2.4);
		\draw (-.5,0) -- (2,0) -- (2,-.5) -- (-.5,2) -- (0,2) -- (0,-.5) -- cycle; 
		\draw [dashed] (-.5,1) -- (1,-.5) -- (1,1) -- cycle;
		\draw (0,-.5) node [left] {$P_{e_1}$};
		\draw (0,2) node [below right] {$P_{e_1-e_3}$};
		\draw (2,0) node [below right] {$P_{e_3}$};
		\draw (-.5,0) node {$\diamondsuit$};
		\draw (-.5,2) node {$\diamondsuit$};
		\draw (0,-.5) node {$\square$};
		\draw (2,-.5) node {$\square$};
		\draw (2,0) node {$\bigcirc$};
		\draw (0,2) node {$\bigcirc$};
		\draw (2,1.1) node {$PP_{\cF}$};
	\end{tikzpicture}\vspace{-8pt}$$
\end{ex}

\begin{ex}
	The following more simple example was pointed out to me by Bivas Khan. Consider the tangent bundle $\cE=\cT_{\cH_2}$ of the Hirzebruch surface, the rays of the fan are $v_1=(1,0)$, $v_2=(0,1)$, $v_3=(-1,1)$, $v_4=(0,-1)$. $$
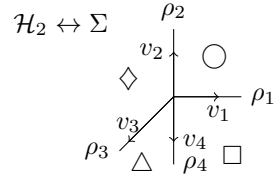
\begin{tikzpicture}[scale=.6]
		\draw (-2.5,1.6) node {$\cH_2\leftrightarrow\Sigma$};
		\draw (0,0) -- (1.5,0) node[right] {$\rho_1$};
		\draw (0,0) -- (0,1.5)node[above] {$\rho_2$};
		\draw (0,0) -- (-1.2,-1.2) node[left] {$\rho_3$};
		\draw (0,0) -- (0,-1.5)node[right] {$\rho_4$};
		\draw [->] (0,0) -- (1,0) node[below] {$v_1$};
		\draw [->] (0,0) -- (0,1) node[left] {$v_2$};
		\draw [->] (0,0) -- (-1,-1) node[above] {$v_3$};
		\draw [->] (0,0) -- (0,-1) node[right] {$v_4$};
		\draw (-1,.4) node {$\diamondsuit$};
		\draw (1.3,-1.3) node {$\square$};
		\draw (.9,.9) node {$\bigcirc$};
		\draw (-.7,-1.5) node {$\triangle$};
	\end{tikzpicture}$$
	The linear space $F=\Span(v_2)$ defines a subbundle $\cF$ and $\cE$ cannot be written as $\cF\oplus\cG$. Here we will see that the obstruction to being a direct sum comes only from the matroid (and not from the form of the polytopes).
	For the picture we tensorize $\cE$ and $\cF$ by the line bundle $\cO(D_3+D_4)$ which does not change the stability.
	$$
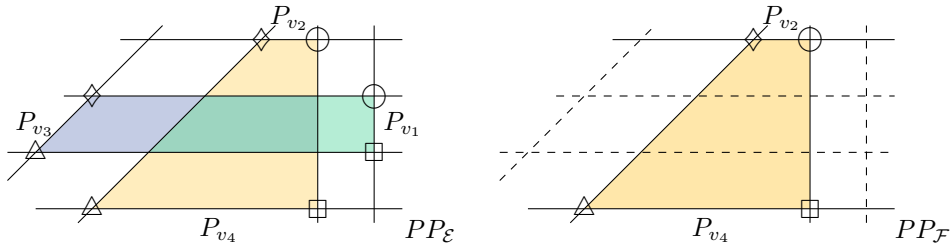
\begin{tikzpicture}[scale=.75]
		\fill [jaune!50, opacity=0.6] (0,-2) -- (-4,-2) -- (-1,1) -- (0,1) -- cycle;
		\fill [bleu!50, opacity=0.6] (0,0) -- (-4,0) -- (-5,-1) -- (0,-1) -- cycle; 
		\fill [vert!50, opacity=0.6] (-2,0) -- (1,0) -- (1,-1) -- (-3,-1) -- cycle;
		\draw (-5,-2) -- (1.5,-2);
		\draw (-5.5,-1) -- (1.5,-1);
		\draw (-4.5,0) -- (1.5,0);
		\draw (-3.5,1) -- (1.5,1);
		\draw (0,-2.25) -- (0,1.25);
		\draw (1,-2.25) -- (1,1.25);
		\draw (-5.5,-1.5) -- (-2.75,1.25);
		\draw (-4.25,-2.25) -- (-.75,1.25);
		\draw (1,-.5) node[right] {$P_{v_1}$};
		\draw (-4.5,-.5) node[left] {$P_{v_3}$};
		\draw (-.5,1) node[above] {$P_{v_2}$};
		\draw (-1.75,-2) node[below] {$P_{v_4}$};
		\draw (-4,0) node {$\diamondsuit$};
		\draw (1,-1) node {$\square$};
		\draw (1,0) node {$\bigcirc$};
		\draw (-4,-2) node {$\triangle$};
		\draw (-1,1) node {$\diamondsuit$};
		\draw (0,-2) node {$\square$};
		\draw (0,1) node {$\bigcirc$};
		\draw (-5,-1) node {$\triangle$};
		\draw (2,-2.4) node {$PP_{\cE}$};
	\end{tikzpicture}\quad \; \begin{tikzpicture}[scale=.75]
		\fill [jaune!40] (0,-2) -- (-4,-2) -- (-1,1) -- (0,1) -- cycle;
		\draw (-5,-2) -- (1.5,-2);
		\draw [dashed] (-5.5,-1) -- (1.5,-1);
		\draw [dashed] (-4.5,0) -- (1.5,0);
		\draw (-3.5,1) -- (1.5,1);
		\draw (0,-2.25) -- (0,1.25);
		\draw [dashed] (1,-2.25) -- (1,1.25);
		\draw [dashed] (-5.5,-1.5) -- (-2.75,1.25);
		\draw (-4.25,-2.25) -- (-.75,1.25);
		\draw (-.5,1) node[above] {$P_{v_2}$};
		\draw (-1.75,-2) node[below] {$P_{v_4}$};
		\draw (-4,-2) node {$\triangle$};
		\draw (-1,1) node {$\diamondsuit$};
		\draw (0,-2) node {$\square$};
		\draw (0,1) node {$\bigcirc$};
		\draw (2,-2.4) node {$PP_{\cF}$};
	\end{tikzpicture}$$
	Remark that we could have taken the same filtrations on any complete smooth toric surface having a fan with $4$ rays for instance $X=\bP^1\times\bP^1$.
\end{ex}

In particular, these examples give rise to nontrivial Harder--Narasimhan filtrations
$$0\subseteq\cF\subseteq\cE\ .$$

\section{Stability of restrictions to invariant curves}\label{s:rest}

The $\alpha$-semistability of the restrictions $\cE|_C$ of a toric vector bundle $\cE$ to a torus invariant curve is much easier to view on the parliament of $\cE$. 
Moreover, it is often useful to compare the $\alpha$-semistability of a toric vector bundle $\cE$ to the semistability of its restrictions $\cE|_C$ to any torus invariant curve. 

\begin{rem}
	Theorem 2.5 from \cite{BLg} furnishes a sufficient condition to deduce from the semistability of a toric vector bundle $\cE$, the semistability of its restriction $\cE|_C$ to any invariant curves $C$ : the characteristic class of $\cE$ is $0$
	$$\Delta(\cE)=c_2(\cE)-\frac{r-1}{2r}c_1(\cE)^2=0\ .$$
\end{rem}

In this section, we explain how to see the (semi)stability of the restrictions $\cE|_C$ of a vector bundle $\cE$ to a torus invariant curve $C$, given a parliament $PP_\cE$. In \cite{DrJS} Subsection 3.1, Di Rocco, Jabbusch and Smith recover the parliaments of restrictions $\cE|_C$ to torus invariant curves from the parliaments of polytopes of $\cE$.

\begin{prop}\label{restr}
	By the cone-orbit correspondence, torus invariant curves correspond to a cone $\tau\in\Sigma(d-1)$. Since $X$ is complete, there are two maximal cones $\sigma$ and $\sigma'$ in $\Sigma(d)$ containing $\tau$. 
	
	If $PP_\cE$ is a parliament of $\cE$, then a parliament of polytopes of $\cE|_C$ is composed of the projection on $\tau^\perp$ of the line segments parallel to $\tau^\perp$ joining the associated characters in $\textbf u(\sigma)$ and $\textbf u(\sigma)$ (renormalized by $1/||u_\tau||$).
\end{prop}

\begin{ex}
	Consider $X=\textup{Bl}_{[0:1:0]}(\bP^2)$ and let $\cE$ be the equivariant vector bundle $$\cE=\cO_X(4D_0+D_3)\oplus\cO_X(3D_1-D_3)$$ on $X$.
	The unique parliament of polytopes of $\cE$ is 
	$$
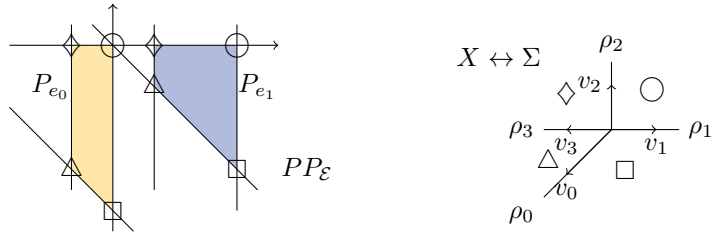
\begin{tikzpicture}[scale=.55]
		\fill [jaune!40] (-1,0) -- (0,0) -- (0,-4) -- (-1,-3);
		\fill [bleu!40] (1,0)-- (3,0) -- (3,-3) -- (1,-1);  
		\draw [->] (-2.5,0) -- (4,0);
		\draw (3,-4) -- (3,.5);
		\draw [->] (0,-4.5) -- (0,1);
		\draw (-.5,.5) -- (3.5,-3.5);
		\draw (-2.5,-1.5) -- (.5,-4.5); 
		\draw (-1,-3.5) -- (-1,.5);
		\draw (1,-3.5) -- (1,.5);
		\draw (3.5,-1) node {$P_{e_1}$};
		\draw (-1.5,-1) node {$P_{e_0}$};
		\draw (4.7,-3) node {$PP_\cE$};
		\draw (-1,0) node {$\diamondsuit$};
		\draw (1,0) node {$\diamondsuit$};
		\draw (3,-3) node {$\square$};
		\draw (0,-4) node {$\square$};
		\draw (3,0) node {$\bigcirc$};
		\draw (0,0) node {$\bigcirc$};
		\draw (-1,-3) node {$\triangle$};
		\draw (1,-1) node {$\triangle$};
	\end{tikzpicture}\qquad \qquad
	\begin{tikzpicture}[scale=.6]
		\draw (-2.5,1.6) node {$X\leftrightarrow\Sigma$};
		\draw (-1.5,0) node[left] {$\rho_3$} -- (1.5,0) node[right] {$\rho_1$};
		\draw (0,0) -- (0,1.5)node[above] {$\rho_2$};
		\draw (-1.5,-1.5) node [below left]  {$\rho_0$}-- (0,0);
		\draw [->] (0,0) -- (1,0) node[below] {$v_1$};
		\draw [->] (0,0) -- (0,1) node[left] {$v_2$};
		\draw [->] (0,0) -- (-1,0) node[below] {$v_3$};
		\draw [<-] (-1,-1) node [below]  {$v_0$} -- (0,0);
		\draw (-1,.8) node {$\diamondsuit$};
		\draw (-1.4,-.7) node {$\triangle$};
		\draw (.3,-.9) node {$\square$};
		\draw (.9,.9) node {$\bigcirc$};
	\end{tikzpicture}$$
	The parliament of $\cE$ restricted to the torus invariant curve $D_0$ is obtained the following way.
	$$
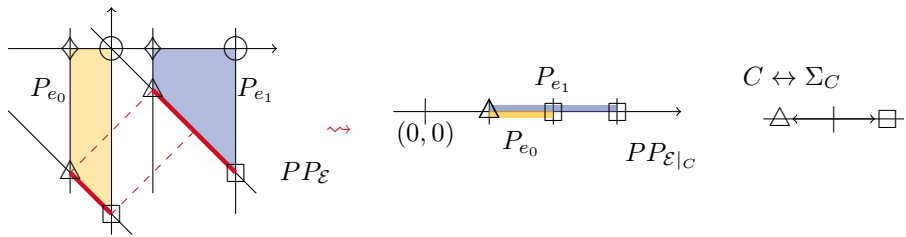
\begin{tikzpicture}[scale=.55]
		\fill [jaune!40] (-1,0) -- (0,0) -- (0,-4) -- (-1,-3);
		\fill [bleu!40] (1,0)-- (3,0) -- (3,-3) -- (1,-1);
		\draw [->] (-2.5,0) -- (4,0);
		\draw (3,-4) -- (3,.5);
		\draw [->] (0,-4.5) -- (0,1);
		\draw (-.5,.5) -- (3.5,-3.5);
		\draw (-2.5,-1.5) -- (.5,-4.5); 
		\draw (-1,-3.5) -- (-1,.5);
		\draw (1,-3.5) -- (1,.5);
		\draw [rouge, line width=1.7] (-1,-3) -- (0,-4);  
		\draw [rouge, line width=1.7] (1,-1) -- (3,-3);  
		\draw [rouge, dashed] (1,-1) -- (-1,-3);  
		\draw [rouge, dashed] (2,-2) -- (0,-4);  
		\draw (3.5,-1) node {$P_{e_1}$};
		\draw (-1.5,-1) node {$P_{e_0}$};
		\draw (4.7,-3) node {$PP_\cE$};
		\draw (-1,0) node {$\diamondsuit$};
		\draw (1,0) node {$\diamondsuit$};
		\draw (3,-3) node {$\square$};
		\draw (0,-4) node {$\square$};
		\draw (3,0) node {$\bigcirc$};
		\draw (0,0) node {$\bigcirc$};
		\draw (-1,-3) node {$\triangle$};
		\draw (1,-1) node {$\triangle$};
		\draw [rouge] (5.5,-2) node {$\rightsquigarrow$};
	\end{tikzpicture}\quad
	\begin{tikzpicture}[scale=.85]
		\draw [->] (-.5,0) -- (4,0);
		\draw (0,0) node {$|$} node [below] {$(0,0)$};
		\draw (1,0) node {$|$};;
		\draw (2,0) node {$|$};
		\draw (3,0) node {$|$};
		\fill [bleu, opacity=0.5] (1,.1) -- (3,.1) -- (3,0) -- (1,0);
		\fill [jaune, opacity=0.7] (1,-.1) -- (2,-.1)-- (2,0) -- (1,0);	
		\draw (3.7,-.7) node {$PP_{\cE|_C}$};
		\draw (0,-1.8) node { };	
		\draw (2,.5) node {$P_{e_1}$};
		\draw (1.5,-.5) node {$P_{e_0}$};
		\draw (1,0) node {$\triangle$};
		\draw (3,0) node {$\square$};
		\draw (1,0) node {$\triangle$};
		\draw (2,0) node {$\square$};
	\end{tikzpicture}\quad
	\begin{tikzpicture}[scale=.55]
		\draw (-1,1) node {$C\leftrightarrow\Sigma_C$};
		\draw (-1.7,0) -- (1.7,0);
		\draw [->] (0,0) node {$|$} -- (1,0);
		\draw [->] (0,0) -- (-1,0);
		\draw (-1.3,0) node {$\triangle$};
		\draw (1.3,0) node {$\square$};
		\draw (0,-2.6) node { };
	\end{tikzpicture}
	$$ \vspace{-13pt}
\end{ex}

The semistability of such restriction $\cE|_C$ should be easier to understand as torus invariant curves $C$ are isomorphic to $\bP^1$. Moreover, in 1957, Grothendieck classified vector bundles on $\bP^1$.
\begin{thm}[\cite{G} Theorem 2.1]
	Every holomorphic vector bundle $\cE$ on $\bP^1$ is holomorphically isomorphic to a direct sum of line bundles:
	$$\cE\cong \cO_{\bP^1}(d_1)\oplus\cdots\oplus\cO_{\bP^1}(d_r)\ .$$
\end{thm}  
The problem of semistability is then reduced to verify if $d_1=...=d_r$. This can be achieved by looking at the parliaments of polytopes of $\cE$. 
The operation on parliaments corresponding to restricting the toric bundle to a invariant curve is explained in Proposition \ref{restr}. We obtain the following result.

\begin{thm}
	Let $C$ be a torus invariant curve on $X$. $C$ corresponds to a wall $\tau\in\Sigma(d-1)$ between two maximal cones $\sigma,\sigma'\in\Sigma(d)$. 
	The restriction $\cE|_C$ of a toric vector bundle $\cE$ to $C$ is semistable if and only if the normalized lattice distance between associated characters $u\in\textbf{u}(\sigma)$ and $u'\in\textbf{u}(\sigma')$ in the one-dimensional lattice $(\tau^\perp+u)\cap M$ is the same for any pair $u,u'$.
\end{thm}

The example which follows concerns a semistable vector bundle such that its restriction to an invariant curve is not semistable.

\begin{ex}
	Consider $X=\textup{Bl}_{[0:1:0]}(\bP^2)$ and let $\cE$ be the equivariant vector bundle $$\cE=\cO_X(4D_0+D_3)\oplus\cO_X(3D_1-D_3)$$ on $X$.
	The parliament $\cE_{D_0}$ of $\cE$ restricted to the torus invariant curve $D_0$ is obtained the following way.
	$$
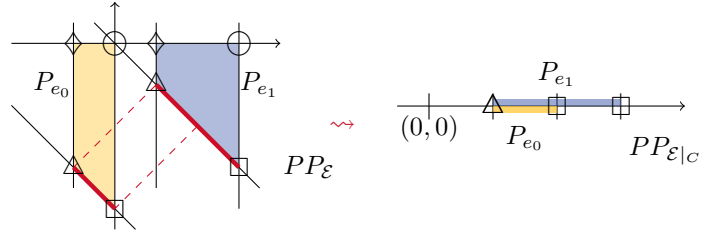
\begin{tikzpicture}[scale=.55]
		\fill [jaune!40] (-1,0) -- (0,0) -- (0,-4) -- (-1,-3);
		\fill [bleu!40] (1,0)-- (3,0) -- (3,-3) -- (1,-1);
		\draw [->] (-2.5,0) -- (4,0);
		\draw (3,-4) -- (3,.5);
		\draw [->] (0,-4.5) -- (0,1);
		\draw (-.5,.5) -- (3.5,-3.5);
		\draw (-2.5,-1.5) -- (.5,-4.5); 
		\draw (-1,-3.5) -- (-1,.5);
		\draw (1,-3.5) -- (1,.5);
		\draw [rouge, line width=1.7] (-1,-3) -- (0,-4);  
		\draw [rouge, line width=1.7] (1,-1) -- (3,-3);  
		\draw [rouge, dashed] (1,-1) -- (-1,-3);  
		\draw [rouge, dashed] (2,-2) -- (0,-4);  
		\draw (3.5,-1) node {$P_{e_1}$};
		\draw (-1.5,-1) node {$P_{e_0}$};
		\draw (4.7,-3) node {$PP_\cE$};
		\draw (-1,0) node {$\diamondsuit$};
		\draw (1,0) node {$\diamondsuit$};
		\draw (3,-3) node {$\square$};
		\draw (0,-4) node {$\square$};
		\draw (3,0) node {$\bigcirc$};
		\draw (0,0) node {$\bigcirc$};
		\draw (-1,-3) node {$\triangle$};
		\draw (1,-1) node {$\triangle$};
		\draw [rouge] (5.5,-2) node {$\rightsquigarrow$};
	\end{tikzpicture}\quad
	\begin{tikzpicture}[scale=.85]
		\draw [->] (-.5,0) -- (4,0);
		\draw (0,0) node {$|$} node [below] {$(0,0)$};
		\draw (1,0) node {$|$};;
		\draw (2,0) node {$|$};
		\draw (3,0) node {$|$};
		\fill [bleu, opacity=0.5] (1,.1) -- (3,.1) -- (3,0) -- (1,0);
		\fill [jaune, opacity=0.7] (1,-.1) -- (2,-.1)-- (2,0) -- (1,0);	
		\draw (3.7,-.7) node {$PP_{\cE|_C}$};
		\draw (0,-1.8) node { };	
		\draw (2,.5) node {$P_{e_1}$};
		\draw (1.5,-.5) node {$P_{e_0}$};
		\draw (1,0) node {$\triangle$};
		\draw (3,0) node {$\square$};
		\draw (1,0) node {$\triangle$};
		\draw (2,0) node {$\square$};
	\end{tikzpicture}
	$$ 
	
	The bundle $\cE|_{D_0}$ would be stable if the segments $P_{e_0}$ and $P_{e_1}$ were of the same length. In this example, there is not any torus invariant curve $C$ such that $\cE|_C$ is semistable.
	A contrario, taking the movable curve $\alpha=2D_1-D_3$, as in Example \ref{blowup}, with its corresponding $(t_i)_{\rho_i\in\Sigma(1)}$ $$
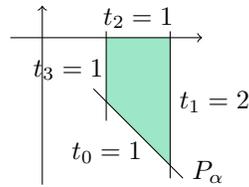
\begin{tikzpicture}[scale=.85]
		\fill [vert!40] (1,0) -- (2,0) -- (2,-2) -- (1,-1); 
		\draw [->] (-.5,0) -- (2.5,0);
		\draw [<-] (0,.5) -- (0,-2.3);
		\draw (.8,-.8) -- (2.2,-2.2);
		\draw (2,.1) -- (2,-2.2);
		\draw (1,.1) -- (1,-1.3);
		\draw (2.2,-2.1) node [right] {$P_\alpha$};
		\draw (1.5,0) node [above] {$t_2=1$};
		\draw (2,-1) node [right] {$t_1=2$};
		\draw (1.1,-.5) node [left] {$t_3=1$};
		\draw (1.7,-1.8) node [left] {$t_0=1$};
	\end{tikzpicture}$$ gives $\cE$ is $\alpha$-polystable: $\mu(\cL_0)=4t_0+t_3=5$ is equal to $\mu(\cL_1)=3t_1-t_3=5$. 
\end{ex}

\begin{ex}
	For any choice of movable curve $\alpha$, the tangent bundle $\cT_{\bP^2}$ of $\bP^2$ is an example of $\alpha$-stable equivariant vector bundle $\cE$ with a non stable restriction to a torus invariant curve.
\end{ex}

\section{Alternative definition of parliaments}\label{s:A}

We propose definitions for the parliaments of polytopes of a toric vector bundle $\cE$. For each definition, the data of a parliament of polytopes of some globally generated equivariant vector bundle $\cE$ corresponds to some isomorphism class of $\cE$ (corresponding to the title of the subsection). 
\begin{rem}
	If we consider virtual polytopes as proposed in Remark \ref{virtual}, the condition of global generatedness becomes superfluous.
\end{rem}

The naive definition for a parliament of polytopes of a toric vector bundle $\cE$ would be
$$PP_\cE:= \left\{P_e\ \middle|\ e\in G(\cE)\right\}$$
which does not keep trace of the label of each polytope. This definition does not even allow us to distinguish a rank-$2$ from a rank-$3$ vector bundle on $\bP^2$ : the tangent bundle $\cT_{\bP^2}$ and the splitting bundle $\cO_{\bP^2}(D_0)\oplus\cO_{\bP^2}(D_1)\oplus\cO_{\bP^2}(D_2)$ have the same parliament of polytopes.

\subsection{Equivariant isomorphism class of framed toric vector bundle}

The definition we gave in Definition \ref{def} for a parliament of polytopes of a toric vector bundle $\cE$ was $PP_\cE:=\left\{(P_e,e)\ \middle|\ e\in G(\cE)\right\}$ where $G(\cE)$ is defined modulo isomorphism of type $(\star)$.

\begin{defn}
	A framed equivariant vector bundle is a toric vector bundle with a choice of isomorphism $E\cong \bC^r$. 
	
	A morphism of framed equivariant vector bundles is a morphism of equivariant vector bundle compatible with the framing.
\end{defn} 

\begin{prop}\label{isom}
	The data of a parliament of polytopes
	$$PP_\cE:=\left\{(P_e,e)\ \middle|\ e\in G(\cE)\right\}\ /(\star)$$
	of a globally generated equivariant vector bundle $\cE$ is equivalent to knowing the equivariant isomorphism class of the framed equivariant vector bundle $\cE$.
\end{prop}

\begin{prr}
	In \cite{P} Proposition 3.4, the equivariant class of a framed equivariant vector bundle is uniquely determined by the following 
	$$\left(\{\textbf u(\sigma)\}_{\sigma\in\Sigma(d)} \ ,\quad \{Fl(\rho)\}_{\rho\in\Sigma(1)}\right)\,,$$
	where $Fl(\rho_i)$ is the flag appearing in the filtration $(E^i(j))_j$.
	From the parliament we recover the flag $$Fl(\rho_i): \{0\}\subset E^i(A^i_r)\subset \ldots \subset E^i(A^i_1)=E$$ by considering all different vector spaces
	$$E_{i,j}=\sum_{\substack{\langle u,v_i\rangle\geq j,\\ u\in P_e \text{ for some }e\in G(\cE)}}e\ .$$
	Indeed, if $e\in E_{i,j}$ then there exists $u\in P_e$ such that $\langle u,v_i\rangle\geq j$ and by definition of $P_e$, we obtain that $e\in E^i(j)$.
	And conversely, by Theorem 1.2 of \cite{DrJS}, as $\cE$ is taken globally generated, any $u\in \textbf u(\sigma)$ belongs to some polytope and by (\ref{CC}), we have $E^i(j)\subseteq E_{i,j}$.
\end{prr}

\subsection{Equivariant isomorphism class of toric vector bundle}

We may not want to deal with framings anymore. Let us thus quotient by the action of $GL_r(\bC)$ on $\bC^r$. 

\begin{prop}
	The data of a parliament of polytopes
	$$PP_\cE:= \left\{(P_e,e)\ \middle|\ e\in G(\cE)\right\}\ /(\star)\ /GL_r(\bC)$$ 
	of a globally generated equivariant vector bundle $\cE$ is equivalent to knowing the equivariant isomorphism class of the equivariant vector bundle $\cE$.
\end{prop}

\begin{rem}
	Here morphisms do not have to preserve the framings anymore.
\end{rem}

\begin{propr}
	It follows from \cite{P} Corollary 3.6 that the equivariant class of a toric vector bundle is uniquely determined by
	$$\left(\{\textbf u(\sigma)\}_{\sigma\in\Sigma(d)} \ ,\quad \cO_{\{Fl(\rho)\}_{\rho\in\Sigma(1)}}\right)\,,$$
	where $\cO_{\{Fl(\rho)\}_{\rho\in\Sigma(1)}}$ is the $GL_r(\bC)$-orbit of the flag given by the filtration $(E^i(j))_j$. We conclude by Proposition \ref{isom}.
\end{propr}

\subsection{Isomorphism class of toric vector bundle}

Finally, we may want an object which represents the isomorphism class of $\cE$ and not its equivariant isomorphism class. 
For that we need to quotient by the group $T$ of compositions of translations for each direct component.

\begin{prop}\label{isom2}
	The data of a parliament of polytopes 
	$$PP_\cE:= \left\{(P_e,e)\ \middle|\ e\in G(\cE)\right\}\ /(\star)\ /GL_r(\bC)\ /T\ ,$$
	of a globally generated equivariant vector bundle $\cE$ is equivalent to knowing the isomorphism class of $\cE$.
\end{prop}

\begin{prr}
	It is known that any line bundle $\cO_X(D)$ on toric variety is isomorphic to an equivariant line bundle $\cO_X(D_T)=\cO_X(D)\otimes\chi^u$ for a unique $u\in M$. 
	
	In fact this result has been generalized by Klyachko in (\cite{K} Corollary 1.2.4). If two equivariant vector bundles $\cE$ and $\cF$ are isomorphic then there exist characters $\chi_1,\hdots,\chi_m$ such that
	\begin{center} $\cE_i\otimes\chi_i$ and $\cF_i$ are equivariantly isomorphic,\ \end{center}
	where $\cE=\cE_1\oplus\ldots\oplus\cE_m$ and $\cF=\cF_1\oplus\ldots\oplus\cF_m$ are some direct decompositions of $\cE$ and $\cF$. 
	
	As the level of parliaments, tensoring by a character $\chi^u$ corresponds to translating the parliament by $u\in M$.
\end{prr}

\newpage
\bibliographystyle{alpha} 
\bibliography{bib}

\Addresses

\end{document}